\documentclass[a4paper]{amsart}
\usepackage[totalwidth=17cm,totalheight=22cm]{geometry}
\makeatletter
\newif\iflabel\labelfalse \let\w@label=\label
\def\label{\global\labeltrue\w@label}
\let\w@eqn=\equation
\def\equation{\global\labelfalse\w@eqn}
\def\endequation{\iflabel\@eeq\else\@eeqw\fi$$\global\@ignoretrue}
\def\@eeq{\eqno \@eqnnum}
\def\@eeqw{\addtocounter{equation}{-1}}
\newskip\bw@\newskip\bws@
\def\dc@@pq{\global \bw@=\belowdisplayskip \global\bws@=\belowdisplayshortskip
   \global\belowdisplayshortskip=3pt plus -3pt
   \global\belowdisplayskip=\belowdisplayshortskip
   \hskip 500pt minus 500pt\relax}
\def\dc@pq{\dc@@pq$$}
\def\fc@pq{$$\hskip 500pt minus 500pt\global
   \belowdisplayshortskip=\bws@\global\belowdisplayskip=\bw@}
\def\coupeq{\dc@pq\fc@pq}
\def\coupepas{\par\noindent\begin{minipage}{\textwidth}}
\def\jusquela{\end{minipage}}
\makeatother

\catcode`\@=11\relax
\newwrite\@unused
\def\typeout#1{{\let\protect\string\immediate\write\@unused{#1}}}
\def\@nnil{\@nil}
\def\@empty{}
\def\@psdonoop#1\@@#2#3{}
\def\@psdo#1:=#2\do#3{\edef\@psdotmp{#2}\ifx\@psdotmp\@empty \else
    \expandafter\@psdoloop#2,\@nil,\@nil\@@#1{#3}\fi}
\def\@psdoloop#1,#2,#3\@@#4#5{\def#4{#1}\ifx #4\@nnil \else
       #5\def#4{#2}\ifx #4\@nnil \else#5\@ipsdoloop #3\@@#4{#5}\fi\fi}
\def\@ipsdoloop#1,#2\@@#3#4{\def#3{#1}\ifx #3\@nnil 
       \let\@nextwhile=\@psdonoop \else
      #4\relax\let\@nextwhile=\@ipsdoloop\fi\@nextwhile#2\@@#3{#4}}
\def\@tpsdo#1:=#2\do#3{\xdef\@psdotmp{#2}\ifx\@psdotmp\@empty \else
    \@tpsdoloop#2\@nil\@nil\@@#1{#3}\fi}
\def\@tpsdoloop#1#2\@@#3#4{\def#3{#1}\ifx #3\@nnil 
       \let\@nextwhile=\@psdonoop \else
      #4\relax\let\@nextwhile=\@tpsdoloop\fi\@nextwhile#2\@@#3{#4}}
\def\psdraft{
	\def\@psdraft{0}
}
\def\psfull{
	\def\@psdraft{100}
}
\psfull
\newif\if@prologfile
\newif\if@postlogfile
\newif\if@noisy
\def\pssilent{
	\@noisyfalse
}
\def\psnoisy{
	\@noisytrue
}
\psnoisy
\newif\if@bbllx
\newif\if@bblly
\newif\if@bburx
\newif\if@bbury
\newif\if@height
\newif\if@width
\newif\if@rheight
\newif\if@rwidth
\newif\if@clip
\newif\if@verbose
\def\@p@@sclip#1{\@cliptrue}
\def\@p@@sfile#1{%\typeout{file is #1}
		   \def\@p@sfile{#1}
}
\def\@p@@sfigure#1{\def\@p@sfile{#1}}
\def\@p@@sbbllx#1{
		\@bbllxtrue
		\dimen100=#1
		\edef\@p@sbbllx{\number\dimen100}
}
\def\@p@@sbblly#1{
		\@bbllytrue
		\dimen100=#1
		\edef\@p@sbblly{\number\dimen100}
}
\def\@p@@sbburx#1{
		\@bburxtrue
		\dimen100=#1
		\edef\@p@sbburx{\number\dimen100}
}
\def\@p@@sbbury#1{
		\@bburytrue
		\dimen100=#1
		\edef\@p@sbbury{\number\dimen100}
}
\def\@p@@sheight#1{
		\@heighttrue
		\dimen100=#1
   		\edef\@p@sheight{\number\dimen100}
}
\def\@p@@swidth#1{
		\@widthtrue
		\dimen100=#1
		\edef\@p@swidth{\number\dimen100}
}
\def\@p@@srheight#1{
		\@rheighttrue
		\dimen100=#1
		\edef\@p@srheight{\number\dimen100}
}
\def\@p@@srwidth#1{
		\@rwidthtrue
		\dimen100=#1
		\edef\@p@srwidth{\number\dimen100}
}
\def\@p@@ssilent#1{ 
		\@verbosefalse
}
\def\@p@@sprolog#1{\@prologfiletrue\def\@prologfileval{#1}}
\def\@p@@spostlog#1{\@postlogfiletrue\def\@postlogfileval{#1}}
\def\@cs@name#1{\csname #1\endcsname}
\def\@setparms#1=#2,{\@cs@name{@p@@s#1}{#2}}
\def\ps@init@parms{
		\@bbllxfalse \@bbllyfalse
		\@bburxfalse \@bburyfalse
		\@heightfalse \@widthfalse
		\@rheightfalse \@rwidthfalse
		\def\@p@sbbllx{}\def\@p@sbblly{}
		\def\@p@sbburx{}\def\@p@sbbury{}
		\def\@p@sheight{}\def\@p@swidth{}
		\def\@p@srheight{}\def\@p@srwidth{}
		\def\@p@sfile{}
		\def\@p@scost{10}
		\def\@sc{}
		\@prologfilefalse
		\@postlogfilefalse
		\@clipfalse
		\if@noisy
			\@verbosetrue
		\else
			\@verbosefalse
		\fi
}
\def\parse@ps@parms#1{
	 	\@psdo\@psfiga:=#1\do
		   {\expandafter\@setparms\@psfiga,}}
\newif\ifno@bb
\newif\ifnot@eof
\newread\ps@stream
\def\bb@missing{
	\if@verbose{
		\typeout{psfig: searching \@p@sfile \space  for bounding box}
	}\fi
	\openin\ps@stream=\@p@sfile
	\no@bbtrue
	\not@eoftrue
	\catcode`\%=12
	\loop
		\read\ps@stream to \line@in
		\global\toks200=\expandafter{\line@in}
		\ifeof\ps@stream \not@eoffalse \fi
		\@bbtest{\toks200}
		\if@bbmatch\not@eoffalse\expandafter\bb@cull\the\toks200\fi
	\ifnot@eof \repeat
	\catcode`\%=14
}	
\catcode`\%=12
\newif\if@bbmatch
\def\@bbtest#1{\expandafter\@a@\the#1%%BoundingBox:\@bbtest\@a@}
\long\def\@a@#1%%BoundingBox:#2#3\@a@{\ifx\@bbtest#2\@bbmatchfalse\else\@bbmatchtrue\fi}
\long\def\bb@cull#1 #2 #3 #4 #5 {
	\dimen100=#2 bp\edef\@p@sbbllx{\number\dimen100}
	\dimen100=#3 bp\edef\@p@sbblly{\number\dimen100}
	\dimen100=#4 bp\edef\@p@sbburx{\number\dimen100}
	\dimen100=#5 bp\edef\@p@sbbury{\number\dimen100}
	\no@bbfalse
}
\catcode`\%=14
\def\compute@bb{
		\no@bbfalse
		\if@bbllx \else \no@bbtrue \fi
		\if@bblly \else \no@bbtrue \fi
		\if@bburx \else \no@bbtrue \fi
		\if@bbury \else \no@bbtrue \fi
		\ifno@bb \bb@missing \fi
		\ifno@bb \typeout{FATAL ERROR: no bb supplied or found}
			\no-bb-error
		\fi
		\count203=\@p@sbburx
		\count204=\@p@sbbury
		\advance\count203 by -\@p@sbbllx
		\advance\count204 by -\@p@sbblly
		\edef\@bbw{\number\count203}
		\edef\@bbh{\number\count204}
}
\def\in@hundreds#1#2#3{\count240=#2 \count241=#3
		     \count100=\count240	% 100 is first digit #2/#3
		     \divide\count100 by \count241
		     \count101=\count100
		     \multiply\count101 by \count241
		     \advance\count240 by -\count101
		     \multiply\count240 by 10
		     \count101=\count240	%101 is second digit of #2/#3
		     \divide\count101 by \count241
		     \count102=\count101
		     \multiply\count102 by \count241
		     \advance\count240 by -\count102
		     \multiply\count240 by 10
		     \count102=\count240	% 102 is the third digit
		     \divide\count102 by \count241
		     \count200=#1\count205=0
		     \count201=\count200
			\multiply\count201 by \count100
		 	\advance\count205 by \count201
		     \count201=\count200
			\divide\count201 by 10
			\multiply\count201 by \count101
			\advance\count205 by \count201
		     \count201=\count200
			\divide\count201 by 100
			\multiply\count201 by \count102
			\advance\count205 by \count201
		     \edef\@result{\number\count205}
}
\def\compute@wfromh{
		\in@hundreds{\@p@sheight}{\@bbw}{\@bbh}
		\edef\@p@swidth{\@result}
}
\def\compute@hfromw{
		\in@hundreds{\@p@swidth}{\@bbh}{\@bbw}
		\edef\@p@sheight{\@result}
}
\def\compute@handw{
		\if@height 
			\if@width
			\else
				\compute@wfromh
			\fi
		\else 
			\if@width
				\compute@hfromw
			\else
				\edef\@p@sheight{\@bbh}
				\edef\@p@swidth{\@bbw}
			\fi
		\fi
}
\def\compute@resv{
		\if@rheight \else \edef\@p@srheight{\@p@sheight} \fi
		\if@rwidth \else \edef\@p@srwidth{\@p@swidth} \fi
}
\def\compute@sizes{
	\compute@bb
	\compute@handw
	\compute@resv
}
\def\psfig#1{\vbox {
	\ps@init@parms
	\parse@ps@parms{#1}
	\compute@sizes
	\ifnum\@p@scost<\@psdraft{
		\if@verbose{
			\typeout{psfig: including \@p@sfile \space }
		}\fi
		\special{ps::[begin] 	\@p@swidth \space \@p@sheight \space
				\@p@sbbllx \space \@p@sbblly \space
				\@p@sbburx \space \@p@sbbury \space
				startTexFig \space }
		\if@clip{
			\if@verbose{
				\typeout{(clip)}
			}\fi
			\special{ps:: doclip \space }
		}\fi
		\if@prologfile
		    \special{ps: plotfile \@prologfileval \space } \fi
		\special{ps: plotfile \@p@sfile \space }
		\if@postlogfile
		    \special{ps: plotfile \@postlogfileval \space } \fi
		\special{ps::[end] endTexFig \space }
		\vbox to \@p@srheight true sp{
			\hbox to \@p@srwidth true sp{
				\hss
			}
		\vss
		}
	}\else{
		\vbox to \@p@srheight true sp{
		\vss
			\hbox to \@p@srwidth true sp{
				\hss
				\if@verbose{
					\@p@sfile
				}\fi
				\hss
			}
		\vss
		}
	}\fi
}}
\catcode`\@=12\relax

\newcommand{\vs}{\vspace{0.3cm}}
\newcommand{\dr}{\partial}
\newcommand{\C}{{\mathbb C}}
\newcommand{\N}{{\mathbb N}}
\newcommand{\R}{{\mathbb R}}
\newcommand{\Z}{{\mathbb Z}}
\newcommand{\II}{I\hspace{-0.1cm}I}
\newcommand{\III}{I\hspace{-0.1cm}I\hspace{-0.1cm}I}
\newcommand{\tr}{\mbox{\rm tr}}
\newcommand{\ric}{\mbox{\rm ric}}
\newcommand{\cotg}{\mbox{\rm cotg}}
\newcommand{\SO}{\mbox{\rm SO}}
\newcommand{\ricb}{\overline{\ric}}
\newcommand{\Rc}{\check{R}}
\newcommand{\Rr}{\stackrel{\circ}{R}}
\newcommand{\diam}{\mbox{\rm diam}}
\newcommand{\area}{\mbox{\rm area}}
\newcommand{\be}{\begin{eqnarray}}
\newcommand{\ee}{\end{eqnarray}}

\def\pointir{\unskip  {. --- \ignorespaces }\hskip0cm}

\newcommand{\deltab}{\overline{\delta}}

\newcommand{\gba}{\overline{g}}
\newcommand{\fb}{\overline{f}}
\newcommand{\hb}{\overline{h}}
\newcommand{\kb}{\overline{k}}
\newcommand{\Db}{\overline{D}}
\newcommand{\Rb}{\overline{R}}
\newcommand{\Kb}{\overline{K}}
\newcommand{\Sb}{\overline{S}}
\newcommand{\nub}{\overline{\nu}}
\newcommand{\Omegab}{\overline{\Omega}}

\newcommand{\pt}{\tilde{p}}
\newcommand{\ut}{\tilde{u}}
\newcommand{\yt}{\tilde{y}}
\newcommand{\Bt}{\tilde{B}}
\newcommand{\Ft}{\tilde{F}}
\newcommand{\Kt}{\tilde{K}}
\newcommand{\ept}{\tilde{\epsilon}}
\newcommand{\gat}{\tilde{\gamma}}
\newcommand{\kappat}{\tilde{\kappa}}
\newcommand{\thetat}{\tilde{\theta}}

\newcommand{\rhot}{\tilde{\rho}}
\newcommand{\taut}{\tilde{\tau}}
\newcommand{\phit}{\tilde{\phi}}
\newcommand{\sit}{\tilde{\sigma}}
\newcommand{\Phit}{\tilde{\Phi}}
\newcommand{\Sigmat}{\tilde{\Sigma}}
\newcommand{\Omt}{\tilde{\Omega}}

\newcommand{\kad}{\stackrel{\bullet}{\kappa}}
\newcommand{\gd}{\stackrel{\bullet}{g}}

\newtheorem{prop}{Proposition}[section]
\newtheorem{df}[prop]{Definition}
\newtheorem{lemma}[prop]{Lemma}
\newtheorem{thm}[prop]{Theorem}
\newtheorem{cor}[prop]{Corollary}
\newtheorem{asser}[prop]{Assertion}
\newtheorem{remark}[prop]{Remark}

\newenvironment{thn}[1]{\vskip 0.2cm \noindent{\bf Theorem #1.} \it}{\rm
\vspace{0.2cm}} 
\newenvironment{lmn}[1]{\vskip 0.2cm \noindent{\bf Lemma #1.} \it}{\rm
\vspace{0.2cm}} 

\newcommand{\btm}{\begin{thm}}
\newcommand{\etm}{\end{thm}}
\newcommand{\blm}{\begin{lemma}}
\newcommand{\elm}{\end{lemma}}
\newcommand{\bcr}{\begin{cor}}
\newcommand{\ecr}{\end{cor}}
\newcommand{\bdf}{\begin{df}}
\newcommand{\edf}{\end{df}}
\newcommand{\bprop}{\begin{prop}}
\newcommand{\eprop}{\end{prop}}
\newcommand{\bas}{\begin{asser}}
\newcommand{\eas}{\end{asser}}
\newcommand{\beq}{\begin{equation}}
\newcommand{\eeq}{\end{equation}}
\newcommand{\bpv}{\begin{proof}}
\newcommand{\epv}{\end{proof}}
\newcommand{\bit}{\begin{itemize}}
\newcommand{\eit}{\end{itemize}}
\newcommand{\bpn}{\begin{pfn}}
\newcommand{\epn}{\end{pfn}}
\newcommand{\btn}{\begin{thn}}
\newcommand{\etn}{\end{thn}}
\newcommand{\bln}{\begin{lmn}}
\newcommand{\eln}{\end{lmn}}

\newenvironment{pfn}[1]{\vskip 0.2cm \noindent{\it Proof #1.}}{$\square$
\vspace{0.2cm}} 

\newcommand{\Omb}{\overline{\Omega}}
\newcommand{\Sib}{\overline{\Sigma}}
\newcommand{\gb}{\overline{g}}
\newcommand{\Ub}{\overline{U}}
\newcommand{\Wb}{\overline{W}}
\newcommand{\db}{\overline{\partial}}

\newcommand{\Met}{\mathcal{M}et}
\newcommand{\Imm}{\mathcal{I}mm}
\newcommand{\CMet}{\mathcal{CM}et}
\newcommand{\cA}{\mathcal{A}}
\newcommand{\cC}{\mathcal{C}}
\newcommand{\cD}{\mathcal{D}}
\newcommand{\cE}{\mathcal{E}}
\newcommand{\cF}{\mathcal{F}}
\newcommand{\cM}{\mathcal{M}}
\newcommand{\cS}{\mathcal{S}}
\newcommand{\cP}{\mathcal{P}}
\newcommand{\cV}{\mathcal{V}}
\newcommand{\cW}{\mathcal{W}}
\newcommand{\CImm}{\mathcal{CI}mm}
\newcommand{\gab}{\overline{\gamma}}
\newcommand{\hyp}{\mathbf{H}^3}
\newcommand{\dhyp}{\partial\hyp}
\newcommand{\cL}{\mathcal{L}}
\newcommand{\isom}{\mathrm{Isom}}
\newcommand{\im}{\mathrm{Im}}
\newcommand{\Na}{\nabla}
\newcommand{\Nat}{\tilde{\nabla}}
\newcommand{\Nab}{\overline{\nabla}}
\newcommand{\Sit}{\tilde{\Sigma}}

\newcommand{\eps}{\epsilon}
\newcommand{\ga}{\gamma}
\newcommand{\si}{\sigma}
\newcommand{\om}{\omega}
\newcommand{\Ga}{\Gamma}
\newcommand{\La}{\Lambda}
\newcommand{\Si}{\Sigma}
\newcommand{\Om}{\Omega}

\begin{document}
\title{Higher Schl{\"a}fli Formulas and Applications II \\
Vector-valued differential relations}
\date{Feb. 2008 (v4)}

\author{Jean-Marc Schlenker}
\address{
Institut de Math{\'e}matiques de Toulouse, UMR CNRS 5219,
Universit{\'e} Toulouse III,
31062 Toulouse Cedex 9, France.}
\email{schlenker@math.ups-tlse.fr; 
http://www.picard.ups-tlse.fr/\~{ }schlenker}

\author{Rabah Souam}
\address{
Institut de Math{\'e}matiques de Jussieu,  CNRS  UMR 7586, Universit{\'e} Paris
7, Case 7012,  2 place Jussieu, 75251 Paris Cedex 05, France.}
\email{souam@math.jussieu.fr}

\begin{abstract}
The classical Schl\"afli formula, and its ``higher''
analogs given in \cite{hsf}, are relations between the variations of the
volumes and ``curvatures'' of faces of different dimensions of a polyhedra
(which can be Euclidean, spherical or hyperbolic) under a first-order
deformation. We describe here analogs of those formulas which are vector-valued
rather than scalar. Some consequences follow, for instance constraints on
where cone singularities can appear when a constant curvature manifold is
deformed among cone-manifolds.
\end{abstract}

\maketitle

%\keywords{polyhedra, Schl{\"a}fli, space-forms, mean curvature}

%\subjclass{52B11}

\section{Introduction and Results}

\subsection{The Schl\"afli formula.}

The celebrated Schl\"afli formula (see e.g. \cite{milnor-schlafli})
relates in a simple way the
variations of the volume and of the dihedral angles, at codimension
2 faces, of a polyhedron under a first-order deformation. In a
$n+1$-dimensional space of constant curvature $K$, the formula reads as:
$$ KdV = \frac{1}{n}\sum_e W_e d\theta_e~, $$
where the sum is over the codimension 2 faces of the polyhedron,
$W_e$ is the volume of the face, and $\theta_e$ is the interior
dihedral angle at that face. 

This simple formula is a key tool in different fields of mathematics,
where it has found and
still finds important applications, for instant in hyperbolic geometry 
(see e.g. \cite{bonahon,boileau-porti,BLP}) or in discrete geometry (see e.g.
\cite{bezdek-connelly2,bezdek-connelly,csikos}). Smooth versions of the Schl\"afli formula
have also appeared as interesting tools for hyperbolic 3-manifolds
or higher dimensional Einstein manifolds (e.g. in 
\cite{sem-era,anderson-L2,volume}).

The Schl\"afli formula is also important in parts of physics, where it is 
known as the ``Regge formula'' and is a key tool in the ``Regge calculus'',
a discretization of gravity based on simplicial metrics (see \cite{regge},  
or \cite{CMS} for some related mathematical questions).

\subsection{The higher Schl\"afli formulas.}

In a previous paper \cite{hsf} we
gave ``higher'' analogs of the Schl\"afli formula:
for each $1\le p \le n-1$,
\begin{equation} \label{eq:E_p}
K(n-p) \sum_j K_j dW_j  + p \sum_i W_i dK_i  = 0~,
\end{equation}
where $j$ runs over the faces of codimension $p$ and
$i$ runs over those of codimension $p+2$,
$W_i$ and $W_j$ denote the volumes of the faces, and $K_i$, $K_j$ their
curvatures.

Our goal here is to give extensions of those formulas which, 
rather than being scalar
equations as the Schl\"afli formula above, are vector-valued. We will also derive 
some consequences. 
The new formulas are not only  extensions of those earlier ones, but
also explain them, since the scalar (higher) Schl\"afli formulas
follow from the vector-valued ones by a simple translation invariance
argument. The vector-valued formulas presented here contain much
more information, however, on the possible variations of the
face areas or curvatures of a polyhedron under a deformation.

\subsection{A simple example.}

We state first an elementary example, which should help
understand the higher-dimensional cases considered below.

\begin{thm} \label{tm:polygon-s}
Let $p$ be a polygon in $S^2$, with vertices $v_1, v_2, \cdots, v_n$, with
edge lengths
$l_1, \cdots, l_n\not\in \{ 0,\pi\}$ and (exterior) angles $\theta_1, \cdots,
\theta_n$. 
For each $i\in \{ 1, \cdots, n\}$, let
$w_i$ be the point in $S^2$ which is dual to the oriented edge $e_i:=(v_i,
v_{i+1})$.
Under any first-order deformation of $p$, the variations of its edge lengths
and  dihedral angles satisfy the equation:
$$ \sum_{i=1}^n v_i d\theta_i - w_i dl_i = 0 ~. $$
Conversely, any first-order variation $(l'_1, \cdots, l'_n, \theta'_1, \cdots,
\theta'_n)$
which satisfies this equation corresponds to a first-order deformation of $p$,
which is
uniquely defined (up to the addition of a trivial deformation).
\end{thm}

The special case of isometric deformations --- for which the terms $l'_i$ are
zero --- has been well-known for some time, see e.g. \cite{alex,gluck-generique}.
The proof, which is given in section 7, follows from Theorem 1.8 below
and from some simple dimension-counting arguments. The only part of this
statement that will be generalized in higher dimensions is the first part,
namely the vector-valued equation which holds during deformations.

It is interesting to note that this formula appears in the  context of classifying 
constant mean curvature (CMC) surfaces in the Euclidean space. It was indeed 
derived independently by K. Grosse-Brauckmann, N. Korevaar, R. Kusner, J. Ratzkin 
and J. Sullivan and used to estimate the dimension of the space of appropriate 
Jacobi fields on these surfaces (\cite{KKR06} uses the formula in  the case of 
isometric first-order deformations and  \cite{GKRS07} uses the general 
case)\footnote{We are grateful to R. Kusner for pointing out those references}. 

Very similar statements hold in the hyperbolic and in the de Sitter space, the
Euclidean 3-dimensional space appearing in Theorem 1.1 is then replaced by the
3-dimensional Minkowski space.
A similar statement also holds for Euclidean polygons in the plane, we state
it here since it is slightly different from the spherical statement.

\begin{thm} \label{tm:polygon-e}
Let $p$ be a polygon in $\R^2$, with vertices $v_1, v_2, \cdots, v_n$, with
edge lengths
$l_1, \cdots, l_n\neq 0$ and (exterior) angles $\theta_1, \cdots, \theta_n$.
For each $i\in \{ 1, \cdots, n\}$, let
$n_i$ be exterior unit normal to the edge $e_i:=(v_i, v_{i+1})$.
Under any first-order deformation of $p$, the variations of its edge lengths
and of its angles satisfy the equations:
\begin{enumerate}
\item $ \sum_{i=1}^n \theta'_i = 0~. $
\item $ \sum_{i=1}^n \theta'_i v_i - l'_i n_i = 0 ~. $
\end{enumerate}
Conversely, any first-order variation $(l'_1, \cdots, l'_n, \theta'_1, \cdots,
\theta'_n)$
which satisfies those equations corresponds to a first-order deformation of
$p$, which is
uniquely defined (up to the addition of a trivial deformation).
\end{thm}

\subsection{Some definitions.}

In the sequel we view $S^{n+1}$ as the 
unit sphere centered at the origin in
$\R^{n+2}$ and we let $x$ denote the position vector in $\R^{n+2}$.

\begin{df}
Let $F$ be a face of $P$. Then:
$$ \pi(F) := \int_F x dv \in \R^{n+2}~, $$
while $V(F)$ is the $n+1-k$-dimensional volume of $F$.
\end{df}

So $\pi(F)/V(F)$ is the barycenter of $F$, considered as a vector of
$\R^{n+2}$.
We also recall here the classical notion of polar duality for spherical
polyhedra.

\begin{df}
Let $P\subset S^{n+1}$ be a (convex) polyhedron, the {\it dual polyhedron}
$P^*$ is the set of points $x\in S^{n+1}$ such that $P$ is contained in the
half-space bounded by
the hyperplane orthogonal to $x$ and not containing $x$. Given a face $F$ of
$P$, the face of $P^*$ dual to $F$ is the set of points $x\in S^{n+1}$ such
that the oriented plane orthogonal to $x$ is a support plane of $P$ at $F$.
\end{df}

It is well-known that $P^*$ is combinatorially dual to $P$ --- faces of
dimension $k$ of $P$ are dual to faces of dimension $n-k$ of $P^*$ --- and
that $(P^*)^*=P$.

Very similar definitions can be used for polyhedra in the hyperbolic space.
Instead of considering $S^{n+1}\subset \R^{n+2}$, one then considers
$H^{n+1}$ as a quadric (a so-called pseudo-sphere) in the $n+2$-dimensional
Minkowski space. The dual polyhedron is  defined as in the spherical
case, but it is contained in de Sitter space $S^{n+1}_1$, a simply connected
(for $n\geq 2$) Lorentz space of constant curvature $1$ defined as
$$ S^{n+1}_1 = \{x\in \R^{n+2}_1 ~ | ~ \langle x,x\rangle =1 \}~. $$
More details can be found in particular in \cite{RH} (or, for smooth
surfaces, in \cite{hmcb}).

\subsection{The first formula.}

We can now state the most general form of the first 
vector-valued Schl\"afli-type
formulas presented here. 

We consider convex polyhedra $P$ in a constant curvature space
$M^{n+1}_K$, with $K\in \{ -1, 1\}$, so that $M^{n+1}_K$ 
is either the hyperbolic 
space or the sphere. 
For each $p\in \{0, \cdots, n+1\}$, we denote by $F_p$ the set of faces of
codimension $p$ of $P$. Similar formulas for the Euclidean space are mentioned 
below. We will use a simple notation.

\begin{df}
Let $P$ be a spherical (resp. hyperbolic) polyhedron, let $F$ be a face of $P$.
We call $\Pi_F$ the orthogonal projection on $vect(F)$, and $\Pi_{F^*}$ the
orthogonal projection on $vect(F^*)$. If $P$ is an Euclidean polyhedron, 
$\Pi_F$ is the orthogonal projection on the codimension $p$ linear space parallel to $F,$ denoted by $F_0,$ and $\Pi_{F^*}$ is the orthogonal
projection on $F_0^\perp$ (the orthogonal complement of $F_0$).
\end{df}

\begin{thm} \label{tm:general}
In any first-order deformation of $P$ and for 
any $p\in \{ 1, 2, \cdots, n-1\}$:
\begin{equation*}
(n-p+1)K\sum_{H\in F_p} V(H^*)\Pi_H(\pi'(H)) 
+ p\left( \sum_{F\in F_{p+2}} V'(F^*) \pi(F) + \sum_{G\in F_{p+1}} V(G)
\Pi_G(\pi'(G^*))\right) = 0~. \tag{\mbox{$E_p$}} 
\end{equation*}
\end{thm}

There is a geometric interpretation of $\Pi_{H}(\pi'(H))$: it is the component of
$\pi'(H)$ corresponding to the deformation of $H$ in the vector subspace 
of $\R^{n+2}$ that it
generates (rather than the component corresponding to the displacement
of $vect(H)$ in $\R^{n+2}$). For any face $F$ of $P$, $V(F^*)$ also has a simple
interpretation: it is the volume of the set of normals to the support planes
of $P$ along $F$, considered as a subset of the dual of the link of $P$ at $F$.

\subsection{The corresponding formulas for $p=0$.}

Formula $(E_p)$ given in Theorem 1.6 takes a simpler form for $p=0$,
where some terms disappear.
 
\begin{thm} \label{tm:codim1}
Let $P\subset M^{n+1}_K$ be a convex polyhedron.
Call $\theta(F)$ the exterior dihedral angle at a face $F$ of codimension 2,
and $G^*\subset S^{n+1}$, if $K=1$(resp. $G^*\subset S^{n+1}_1$, if $K=-1$),
the point dual to the face $G$ of codimension
$1$. Then, under any first-order deformation of $P$: 
\begin{equation*} 
\sum_{F\in F_2} \theta'(F) \pi(F) - \sum_{G\in F_1} V'(G) G^* = 0~.
\tag{\mbox{$E_0$}}
\end{equation*}
\end{thm}

The proof can be found in section 5. 

\subsection{Another type of formula.}

As a consequence of those formulas, we find other vector-valued Schl\"afli
type formulas, which are simpler in that they do not involve the projection
terms in formula $(E_p)$ but, on the other hand, they 
involve faces of four different dimensions rather
than only three in $(E_p)$. 

\begin{thm} \label{tm:quatre}
Let $P\subset M^{n+1}_K$ be a convex polyhedron, and let $p\in \{ 2, \cdots, n-1\}$.
For any first-order deformation of $P$,
\begin{equation*}
(n-p+1)K\left( \sum_{F\in F_{p-1}} V'(F) \pi(F^*) + \sum_{G\in F_p} V(G^*)
 \pi'(G)\right) +
p\left( \sum_{H\in F_{p+1}} V(H) \pi'(H^*) + \sum_{L\in F_{p+2}} 
V'(L^*) \pi(L) \right) = 0~. 
\tag{\mbox{$H_p$}}
\end{equation*}
\end{thm}

There is an apparent contradiction between this formula and formula $(E_p)$ 
above since, for $K=0$, the remaining terms look different. In fact $K=0$
corresponds to the Euclidean case, then formula $(M'_p)$ below shows that
the two terms actually correspond.

Those formulas can be better understood, and formulated in slightly
different ways, using the following simple remark. The proof is given
in section 5. 

\begin{remark} \label{rk:18}
\begin{enumerate}
\item For any polyhedron $P\subset M_K^{n+1}$ oriented by the exterior unit
  normal, we have: 
$$ (n+1)K \int_{int(P)} xdv + \sum_{F\in F_1} V(F) F^* = 0~. $$
\item Suppose $P\subset M^{n+1}_K$ is a convex polyhedron. Then for $p=1, 2
  \cdots, n$, we have: 
$$ (n-p+1)K \sum_{F\in F_{p}} V(F^*) \pi(F) + p \sum_{G\in F_{p+1}} V(G) \pi(G^*)
= 0~, $$
with the convention that for a face $G$ of codimension $n+1$, 
that is, for a vertex of $P$, $V(G)=1$.
\end{enumerate}
\end{remark}

It is now possible to remove from formula $(H_p)$ the terms involving $\pi'$. 

\begin{cor} \label{cr:quatre}
  Let $P\subset M^{n+1}_K$ be a convex polyhedron, and let $p\in \{ 2, \cdots, n-1\}$.
For any first-order deformation of $P$,
\begin{equation*}
(n-p+1)K\left( \sum_{F\in F_{p-1}} V'(F) \pi(F^*) - \sum_{G\in F_p} V'(G^*)
 \pi(G)\right) +
p\left( - \sum_{H\in F_{p+1}} V'(H) \pi(H^*) + \sum_{L\in F_{p+2}} 
V'(L^*) \pi(L) \right) = 0~. 
\tag{\mbox{$K_p$}}
\end{equation*}
\end{cor}

This corollary follows directly from Theorem \ref{tm:quatre} and from the
second point in Remark \ref{rk:18}.

\subsection{The Euclidean space.}

As in many cases one can obtain Euclidean versions of the spherical (or
hyperbolic) results given above by considering Euclidean polyhedra,
scaling them by a factor going to $0$, and sending them on $S^{n+1}$ by the
projective map. Precise arguments can be found in section 5, we state the
results here. 

\begin{thm} \label{tm:codim1-e}
Let $P\subset \R^{n+1}$ be a polyhedron.
Call $\theta(F)$ the exterior dihedral angle at a face $F$ of codimension 2,
and, for any codimension $1$ face $G$ of $P$, let $N(G)\in S^{n+1}$ be the
unit vector orthogonal to $G$ towards the exterior of $P$.
Then, under any first-order deformation of $P$: 
\begin{equation*}
\sum_{F\in F_2} \theta'(F) \pi(F) - \sum_{G\in F_1} V'(G) N(G) = 0~.
\end{equation*}
\end{thm}

There is also a direct analog of Theorem \ref{tm:general}. 
We need to adapt some of the definitions to the Euclidean setting.
Let $P\subset \R^{n+1}$ be a polyhedron, let $F$ be a codimension
$p$ face of $P$, and consider a first-order deformation of $P$ (and thus
of $F$). $F^*$ is defined mostly as in the spherical (or hyperbolic)
context, it is now a subset of the unit sphere $S^n$, actually the
interior of a spherical polyhedron. 

We denote by $F_0$ the codimension $p$ plane parallel to $F$ going 
through the origin (in other terms, $F_0$ is a linear subspace, rather
than an affine subspace) and by $F_\perp$ its orthogonal complement. 

\begin{thm} \label{tm:general-e}
In any first-order deformation of $P$ and for 
any $p\in \{ 1, 2, \cdots, n-1\}$:
$$
\sum_{F\in F_{p+2}}
V'(F^*) \pi(F) + \sum_{G\in F_{p+1}} V(G)\Pi_G(\pi(G^*)) = 0 ~.
$$
\end{thm}

It is not quite evident that this formula is invariant under translation
of $P$ (a condition that is obviously necessary). Actually the right-hand
term has this invariance property by the definition, while
the left-hand term changes, under a translation of vector $a$, by
$$ a\sum_{F\in F_{p+2}} V'(F^*)V(F)~. $$
The fact that this sum is zero is precisely the Euclidean
version of the higher Schl\"afli formula of \cite{hsf}. 

Other examples can be obtained,  in particular analogs of the spherical
formulas in Theorem \ref{tm:quatre} and in Corollary \ref{cr:quatre}, can be obtained by
just tingke $K=0$ in those statements.

\subsection{Non-convex polyhedra}

It is interesting to point out the following difference between the classical Schl{\"a}fli formula and ours  (either scalar or vector-valued). The classical Schl{\"a}fli formula is readily seen to be {\it linear}. As a consequence it needs only  to be proved  for  simplices and extends to arbitrary poyhedra by linearity. Our generalized formulas do not have this linearity property. The reason is that, unlike the dihedral angle at codimension 2 faces, the volume and the   mapping  $\pi$ evaluated on polar duals are not linear. Nevertheless it is possible to extend our formulas to general polyhedra once they are known for simplices. Indeed, consider a  convex polyhedron  $P$ which is cut  into two convex pieces $P'$ and $P''$. Then given a face $F$  common to $P'$ and $P''$, it is possible to relate the volumes 
of its polar duals (resp. the mapping $\pi$ evaluated on  its polar duals) 
in the polyhedra $P, P',P''$ and $P'\cap P''$, the latter being seen as a polyhedron in a hyperplane 
(cf. Proposition \ref{nonconvex}). 

Furthemore this observation allows 
us to extend  our formulas to  non-convex polyhedra. 
We explain this with more details in section 6, some of the definitions have to be made
with some care. This might be interesting in view of possible applications
to questions of rigidity for non-convex polyhedra (in dimension 3 or
higher).

\subsection{Applications to deformations of cone-manifolds.}

It is possible to deduce from the vector-valued 
Schl\"afli formulas described here
some interesting properties of singular deformations of spherical, Euclidean
or hyperbolic manifolds. The statements given here are simple and direct 
applications of the vector-valued Schl\"afli formulas that can be found
above.

We describe first the simplest situation, in which we
deform $S^2$ so that cone singularities appear, to show how the well-known
Kazdan-Warner obstruction appears naturally. 

We consider a one-parameter family of spherical metrics $(g_t)_{t\in [0,1]}$
on $S^2$, with cone singularities at points $v_1, \cdots, v_n$, where the cone
angle is equal to $\theta_i(t)$, and we suppose that $g_0$ is the canonical
metric, i.e., $\theta_i(0)=2\pi$ for all $i\in \{1, \cdots, n\}$.

\begin{thm} \label{tm:cone-s-2}
Under those hypothesis,
\be \label{eq:s2}
  \sum_{i=1}^n \theta'_i(0) v_i=0~,
\ee
where the $v_i$ are considered as points in $S^2\subset \R^3$. 
Moreover, given points $v_i$ and numbers $t_i$ such that $\sum_{i=1}^n t_i
v_i=0$, there is a unique first-order deformation of $S^2$ as a cone-manifold
with cone points at the $v_i$ and first-order variations of the cone angles equal
to the $t_i$.
\end{thm}

The previous theorem is actually also a direct consequence of (a singular
version of) the Kazdan-Warner relation \cite{kazdan-warner}. 
In higher dimensions, 
there is a simple relation involving only faces of codimension $2$.
For a face $F$ of codimension $2$ of a cell decomposition $\cC$ of $S^{n+1}$, it is quite natural to call $K(F)$
the ``curvature'' at $F$, i.e., $2\pi$ minus the sum of the angles at $F$ of
the cells of $\cC$ containing it. This singular curvature is $0$ for the spherical
cone metric obtained from a cell decomposition of the sphere, since the angles
then add up to $2\pi$, but it becomes non-zero as soon as this metric is deformed
so that singularities occur along codimension 2 faces.

\begin{thm} \label{tm:sphere-codim1}
Let $\cC$ be a cell decomposition of $S^{n+1}$, let $C_k$ be the set of faces
of codimension $k$ of $\cC$. Consider a first-order deformation of the canonical
metric on $S^{n+1}$
among spherical cone-manifolds with cone singularities at the faces of $\cC$.  Then:
$$ \sum_{F\in C_2} K'(F) \pi(F) = 0~. $$
\end{thm}

The simplest illustration is obtained for $n=2$, for a deformation of
the sphere for which cone singularities appear along a graph. The 
theorem then states that the singular graph, ``weighted'' by the first-order
variations of the singular curvatures at the edges, has to be ``balanced''.

It is possible to extend the previous statement to higher codimensions, 
based on Theorem \ref{tm:quatre}, yielding a wealth of
relations between the possible variations of the volumes and ``curvatures''
of faces of different dimensions of $\cC$. A simple definition is
needed. 

\begin{df}
Let $\cC$ be a cell decomposition of $S^{n+1}$, i.e., a decomposition of
$S^{n+1}$ as the union of spherical  polyhedra (glued isometrically along
their faces). Let $F$ be a codimension $p$ face of $\cC$, $1\leq p\leq n+1$.
We define:
\begin{itemize}
\item $W^*(F):=\sum_{C\supset F} V(F^*_C)$,
\item $F^\circ:=\sum_{C\supset F} \pi(F^*_C)$,
\end{itemize}
where the sum is over the (maximal dimension) cells of $\cC$ containing 
$F$, and $F^*_C$ is the face dual of $F$, considered as a face of $C$ (so
that $F^*_C$ is a face of $C^*$).
\end{df}

By construction $F^\circ$ is a vector in $\R^{n+2}$ orthogonal
to $vect(F)$.

\begin{thm} \label{tm:sphere}
Under the same hypothesis as for Theorem \ref{tm:sphere-codim1},
and for all $p\in \{2, \cdots, n-1\}$:
$$ 
(n-p+1)\left( \sum_{F\in F_{p-1}} V'(F) F^\circ - \sum_{G\in F_p} W^*(G)'
 \pi(G)\right) +
p\left( - \sum_{H\in F_{p+1}} V'(H) H^\circ + \sum_{L\in F_{p+2}} 
W^*(L)' \pi(L) \right) = 0~. 
$$ 
\end{thm}

The proof, which is given in section 7, follows quite directly from Theorem
\ref{tm:quatre}. It is also possible to use Theorem \ref{tm:general} or
Corollary \ref{cr:quatre} to obtain
other, more or less similar statements, we leave this point to the reader. 

Those statements are not restricted to the spherical setting, however the fact
that the manifold one starts with is simply connected is important. It remains
easy to state similar results for deformations, again among cone-manifolds, of
hyperbolic or Euclidean polyhedra; in that case boundary terms appear, unless
one makes relevant hypothesis to the extent that some boundary quantities are
fixed. 

\begin{thm} \label{tm:polygone-h}
Let $p$ be a convex hyperbolic polygon, with vertices $v_1, \cdots, v_n$ and
edges $e_1, \cdots, e_p$. Let $w_1, \cdots, w_q$ be distinct points in the
interior of $p$, and let $l'_1, \cdots, l'_n$, 
$\alpha'_1, \cdots, \alpha'_p$ and 
$\theta'_1, \cdots, \theta'_q$ be real numbers. There exists a first-order
deformation of $p$ among hyperbolic polygons with cone singularities 
at the $w_i$, with:
\begin{itemize}
\item the first-order variation of the length of $e_i$ equal to $l'_i$, 
$1\leq i\leq n$,
\item the first-order variation of the angle at $v_j$ equal to $\alpha'_j$, 
$1\leq j\leq p$,
\item the first-order variation of the total curvature at $w_k$ equal to 
$\theta'_k$, $1\leq k\leq q$, if and only if 
$$ \sum_{i=1}^n l'_i e_i^* + \sum_{j=1}^p \alpha'_j v_j + 
\sum_{k=1}^q \theta'_k w_k =0~. $$
\end{itemize}
This first-order deformation is then unique.
\end{thm}

This simple statement could be considered as a kind of extension of the
Kazdan-Warner obstruction to conformal deformations of hyperbolic
polygons. The proof is in section 7. It should also be possible to
give  hyperbolic analogs of Theorems \ref{tm:cone-s-2} and 
\ref{tm:sphere} in this context of
deformation of hyperbolic polyhedra with cone singularities appearing
in the interior.

It also appears possible that related statements exist for deformations
of hyperbolic or Euclidean manifolds (or cone-manifolds). However stating
such statements would certainly require some care. We believe that it
would be necessary to consider formulas taking values not in $\R^{n+2}$ 
or $\R^{n+2}_1$, as in the simpler cases considered here where the 
fundamental group is trivial, but rather in a $n+2$-dimensional flat
bundle over the manifold considered (or over the regular part of the
cone-manifold under considerations). We do not go further in this
direction here, since it would take us too far from our original
motivations.

\section{Low-dimensional examples}

The geometric meaning of formulas $(E_p), (E'_p)$ and $(H_p)$ is not
obvious at first sight. In this section we describe in
more details what happens in dimension 3 and dimension 4, in the hope of
offering the reader a better grasp of those formulas. 

\subsection{3-dimensional polyhedra}

Obviously formula
$(H_p)$ does not apply to polyhedra in 3-dimensional spaces, since then
$n=2$ and $p$ has to be between $2$ and $n-1$. So we consider in more
details formulas $(E_1)$ for 3-dimensional polyhedra in this subsection,
then formula $(H_2)$ for 4-dimensional polyhedra in the next subsection.
Formula $(E_0)$ is quite simple anyway so we do not consider it here.

Let $P\subset M^3_K$ be a convex polyhedron. We first note that 
if $v$ is a vertex of $P$, then $V(v^*)$ has a simple interpretation.

\begin{remark}
$V(v^*)$ is the singular curvature at $v$ of the induced metric on the
boundary of $P$.
\end{remark}

\begin{proof}
By definition $V(v^*)$ is the area of the face $v^*$ of $P^*$ dual to $v$.
But, by the definition of the duality, $v^*$ is isometric to the dual
(in $S^2$) of the link of $P$ at $v$. So, by the Gauss formula, its
area is $2\pi$ minus the sum of the angle so the dual of the link of 
$P$ at $v$. The angles of the dual of the link are equal to the lengths
of the edges of the link of $P$ at $v$, that is, to the angles at $v$
of the faces of $P$ containing it. So $V(v^*)$ is equal to $2\pi$
minus the sum of the total angle at $v$, or in other terms to the
singular curvature at $v$ of the induced metric on the boundary of $P$.
\end{proof}

Note that this argument extends to the case where $v$ is a codimension 
$3$ face of $P\subset M^{n+1}_K$, for $n\geq 2$. In this subsection, we 
use the notation $k(v)$ for the singular curvature at $v$ of the 
induced metric on $P$, so that $k(v)=V(v^*)$. 

We now consider what Theorem \ref{tm:general} means, first in the
Euclidean 3-dimensional space. Then $n=2$, and the only possible value of 
$p$ is $p=1$. The first term in formula $(E_1)$ vanishes since $K=0$, 
the second sum is over vertices, and the right-hand side is a sum over
edges. We call $\cV, \cE$ the sets of vertices and edges of $P$, respectively.
For each edge $e$, $\Pi_e$ is a linear map from $e_0$ (the linear
line parallel to $e$) to $e^\perp$, and $(E_1)$ reads as follows.

\begin{cor}
In any first-order deformation of a polyhonedron in $\R^3$,
$$ \sum_{v\in \cV} k'(v) v =- \sum_{e\in \cE} l(e) \Pi_e(\pi'(e^*))~. $$
\end{cor}

The meaning of the first sum is clear, while the second sum involve
the way the edges ``rotate''. The statement is illustrated in Figure 1,
with signs indicating how the curvature at a vertex varies. 

\begin{figure}[ht]
\centerline{\psfig{figure=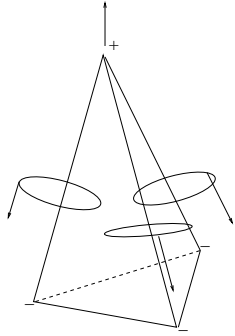,width=3cm}}
\caption{The deformation of a Euclidean simplex.} 
 \label{fig:tetra}
\end{figure}

In the sphere or the hyperbolic space the first term in formula $(E_1)$
kicks in, but in a simpler form since $V(H^*)=1$ when $H$ is codimension
1 face ($H^*$ is then a point). The terms $\pi(H)$ have a very simple
meaning, they are just the barycenter of the faces (in the ambient
4-dimensional space, $\R^4$ or $\R^4_1$) times the area of $H$. 
Theorem \ref{tm:general} then translates as follows.

\begin{cor}
Let $P\subset S^3$ (resp. $H^3$), let $\cF$ be the set of 2-dimensional faces of $P$. 
Then
$$ 2K\sum_{f\in \cF} \Pi_{f}(\pi'(f)) + \sum_{v\in \cV} k'(v) v
=  -\sum_{e\in \cE} l(e) \Pi_e(\pi'(e^*)) ~. $$
\end{cor}

Theorem \ref{tm:quatre} is empty in dimension 3 since $p$ has to be 
chosen between $2$ and $n-1=1$. We do not elaborate on Theorem \ref{tm:codim1-e}
since it is quite clear.

\subsection{Dimension four}

We now turn to $4$-dimensional polyhedra, first in the Euclidean and then
in the spherical (or hyperbolic) setting. 

In $\R^4$, Theorem \ref{tm:general} can be applied either with $p=1$ or with
$p=2$. Now $\cF_2$ the set of $2$-faces of $P$, and, for $f\in \cF_2$,
$a(f)$ is its area. When $e$ is an edge of $P$, we still use the notation
$k(e)=V(e^*)$ since, as pointed out above, $V(e^*)$ is the singular 
curvature at $e$ of the induced metric on the boundary of $P$.

\begin{cor}
Let $P\subset \R^4$, formula $(E_1)$ becomes
$$ \sum_{e\in \cE} k'(e) \pi(e) = -\sum_{f\in \cF_2} a(f) \Pi_f(\pi'(f^*))~, $$
while formula $(E_2)$ can be stated as 
$$ \sum_{v\in \cV} V'(v^*) v = -\sum_{e\in \cE} l(e) \Pi_e(\pi'(e^*))~. $$
\end{cor}

Theorem \ref{tm:quatre} can be applied in $\R^4$ with $p=2$ only, it gives
a formula similar to the one obtained in the previous corollary for $p=2$
but slightly simpler, the two formulas are actually equivalent once one
takes into account formula $(M'_2)$ below.

In $S^4$ or $H^4$, Theorem \ref{tm:general} takes a slightly more complicated
form. For each edge $e$ of $P$, we call $\theta(e)$ the exterior dihedral
angle of $P$ at $e$, and we also call $\cF_3$ the set of $3$-dimensional
faces of $P$.

\begin{cor}
Let $P\subset S^4$ (resp. $P\subset H^4$), formula $(E_1)$ becomes
$$ 3K\sum_{f\in \cF_3} \Pi_f (\pi'(f)) + 
\sum_{e\in \cE} k(e) \pi(e) = -\sum_{f\in \cF_2} a(f) \Pi_f(\pi'(f^*))~, $$
while $(E_2)$ can be written as
$$ 2K\sum_{f\in \cF_2} \theta(f) \Pi_f(\pi'(f)) + 
\sum_{v\in \cV} V'(v^*) v = -\sum_{e\in \cE} l(e) \Pi_e(\pi'(e^*))~. $$
\end{cor}

Finally in this $4$-dimensional context Theorem \ref{tm:quatre} can be
applied non-trivially for $p=2$ to obtain the following result. 

\begin{cor}
Let $P\subset S^4$ (resp. $P\subset H^4$), then, in a first-order deformation, 
$$ K\left( \sum_{f\in \cF_3} V'(f) f^* +\sum_{h \in\cF_2}\theta(h) \pi'(h) \right) + 
\left( \sum_{e\in \cE} l(e) \pi'(e^*) + \sum_{v\in \cV} V'(v^*)v \right) =0~. $$
\end{cor}

\section{Smooth differential formulas}
 
The polyhedral formulas given in the introduction are proved here using analogous
formulas for first-order deformations of smooth hypersurfaces. Let
$\phi:\Sigma\longrightarrow M_K^{n+1}$ be an 
immersed oriented hypersurface in $M_K$, $K\in \{ -1,0,1\}$, 
and let $N$ be an oriented unit
normal vector field along $\Sigma$. Without loss of generality, we will
sometimes, when dealing with local matters, implicitely identify 
the immersion $\phi$ with an inclusion map. Let
$\overline{ \nabla}$ denote the connection on $M_K^{n+1}$.
We define the shape operator of $\Sigma$ as:
$$ Bu =\overline{ \nabla}_uN~, $$
where $u$ is a vector tangent to $\Sigma$. Note that this sign convention is
opposite the one we used in \cite{hsf} and is such that $\langle Bx,x\rangle$
is non-negative for convex 
hypersurfaces if $N$ is oriented towards the exterior (i.e., towards the
concave side of $\Sigma$).

\subsection{Notations}

In all the paper, we call $\nabla$ the Levi-Civit\`a connection of
the hypersurface $\Sigma\subset M_K^{n+1}$, as well as its natural extension
to tensor fields on $\Sigma$. So, given a function $f:\Sigma\rightarrow
\R$, $\nabla f$ is its gradient on $\Sigma$. Given a vector field $Y$
tangent to $\Sigma$, its divergence is defined as:
$$ div(Y) := \sum_{i=1}^n \langle \nabla_{e_i}Y,e_i\rangle~, $$
where $(e_i)_{i=1,\cdots,n}$ is an orthonormal moving frame on $\Sigma$.
Given a vector-valued 1-form $T$ on $\Sigma$, its divergence is
similarly defined as:
$$ div(T)(Z) = \left\langle \sum_{i=1}^n (\nabla_{e_i}T)(Z), e_i
\right\rangle~. $$
We also call $D$ the flat connection on $\R^{n+2}$ and $\R^{n+2}_1$.
Let $X$ be a vector field tangent to
$M_K^{n+1}$ defined on $\Sigma$, which we consider as a first-order deformation
of $\Sigma$. It induces a first-order deformation of the extrinsic invariants
defined on $\Sigma$.

\subsection{The Newton operators.}

We will use the elementary symmetric functions
$S_r$ of the principal curvatures $k_1, k_2,\ldots,k_n$ of the immersion
$\phi$ :
$$
S_r = \sum_{i_1<\ldots<i_r} k_{i_1} \ldots k_{i_r} \,\,\, (1\le r\leq n)~,
$$
as well as the Newton operators, defined as follows.
\begin{df}
The Newton operators of $\Sigma$ are defined for
$0\leq r \leq n$, as:
$$ T_r = S_r Id - S_{r-1}B + \ldots  (-1)^r B^r~, $$
or, inductively, by  $T_0 = Id, T_r = S_r Id - BT_{r-1}$.
\end{df}

\begin{lemma}\label{newton}
The Newton operators satisfy the following formulas for $0\leq r\leq n-1$:
\begin{enumerate}
\item $div(T_r)=0$.
\item $Trace (T_r) = (n-r) S_r $
\item $Trace (BT_r) = (r+1) S_{r+1}$
\item  $S^{\prime}_r= Trace (B^{\prime} T_{r-1})$, the derivative being
taken with respect to any first-order deformation  of the hypersurface.
\end{enumerate}
\end{lemma}

The proofs can be found in \cite{Reilly,rosenberg-constant}.

\subsection{Smooth formulas in non-zero curvature.}

Let $\xi$ be the orthogonal projection of $X$ on $\Sigma$, and let $f:=\langle
X, N\rangle$, where $N$ is the unit vector field normal to $\Sigma$, so that
$X=\xi + fN$. Here we denote by $I$ the induced metric on $\Sigma$, which is 
also classically known as its first fundamental form.

\begin{lemma} \label{lm:23}
For $p=0, \cdots, n-1$:
\begin{equation*}\tag{\mbox{$F_p$}}
\int_\Sigma T_p \xi + \frac{x}{2} \langle I', T_p\rangle -  (p+1) xfS_{p+1}
dv = 0~.
\end{equation*}
\end{lemma}

\begin{proof}
Let $a\in \R^{n+2}$ if $K=1$ (resp. $a\in \R^{n+2}_1$ if $K=-1$)
be a fixed vector. The proof is based on the
integration over $\Sigma$ of the function defined at a point $x\in \Sigma$
as $div(\langle x,a\rangle T_p\xi)$. The vector field $X$ is defined
only on $\Sigma$, but we extend it as a smooth vector field in a
neighborhood of $\Sigma$.

First note that the gradient of $\langle x,a\rangle$ on $\Sigma$ is equal
to the orthogonal projection of $a$ on $\Sigma$, which we call $\overline{a}$.
Indeed, given an orthonormal moving frame $(e_i)_{i=1,\cdots,n}$ on $\Sigma$, 
we have:
$$ \nabla \langle x,a\rangle = \sum_{i=1}^n (e_i.\langle x,a\rangle) e_i
= \sum_{i=1}^n \langle D_{e_i}x, a\rangle e_i
= \sum_{i=1}^n \langle e_i, a\rangle e_i = \overline{a}~. $$
It follows that:
$$ div(\langle x,a\rangle T_p\xi) = \langle \nabla \langle x,a\rangle,
T_p\xi\rangle + \langle x,a\rangle div(T_p\xi) = \langle a, T_p\xi \rangle
+ \langle x,a\rangle div(T_p\xi)~. $$
But:
\begin{eqnarray*}
div(T_p\xi) & = & \sum_{i=1}^n \langle \nabla_{e_i}(T_p\xi), e_i\rangle \\
& = & \sum_{i=1}^n \langle (\nabla_{e_i}T_p)\xi +
T_p \nabla_{e_i}\xi , e_i\rangle \\
& = & \left\langle \sum_{i=1}^n (\nabla_{e_i}T_p)e_i, \xi\right\rangle
+ \langle \nabla_{e_i}\xi , T_p e_i\rangle \\
& = & div(T_p)( \xi)
+\sum_{i=1}^n  \langle  \overline{\nabla}_{e_i}(X-fN) , T_p e_i\rangle \\
& = &\sum_{i=1}^n \langle  \overline{\nabla}_{e_i}X - fBe_i, T_p e_i\rangle \\
& = & \sum_{i=1}^n\langle  \overline{\nabla}_{e_i}X, T_p e_i\rangle - f Trace(BT_p)~. 
\end{eqnarray*}

Let $u,v$ be two vector fields tangent to $\Sigma$. 
Let $\phi_t,\, t\in(-\epsilon,\epsilon),\, \epsilon >0,$ 
be a deformation of the immersion $\phi$ having $X$ as 
a deformation vector field. Then we can write

\begin{eqnarray*}
X\langle u, v\rangle &=& I^{\prime} (u,v) \\
&=&\left( {\frac{d}{dt}}\right)_{|t=0} \langle d\phi_t(u), d\phi_t(v)\rangle\\
&=& \left\langle \left({\frac{D}{dt}}\right)_{|t=0} d\phi_t (u) , 
v\right\rangle + \left\langle u, 
\left({\frac{D}{dt}}\right)_{|t=0} d\phi_t (v)\right\rangle\\
&=& \langle \overline{\nabla}_u X, v\rangle + \langle u,  
\overline{\nabla}_v X\rangle ~. 
\end{eqnarray*}
It follows that the symmetrization of $ \overline{\nabla} X$ is equal to half
the first-order variation of $I$. Therefore:
\be \label{eq:divxi}
 div(\langle x,a\rangle T_p\xi) = \left\langle a, T_p \xi +
\frac{x}{2}\langle
I', T_p\rangle - f(p+1) x S_{p+1}\right\rangle~.
\ee
 The result now follows by integration of this equality over $\Sigma$.
\end{proof}

\begin{lemma} \label{lm:24}
Let $\Sigma\in M_K^{n+1}$, $K\in \{-1,1\}$.
For $p=0, \cdots, n-1$:
\begin{equation*}\tag{\mbox{$G_p$}}
\int_\Sigma \left(
S'_{p+1} + \frac{1}{2} \langle I', BT_p\rangle + (n-p) Kf S_p\right)
x + T_p \overline{\nabla}_XN dv = 0~,
\end{equation*}
where $\overline{\nabla}_XN$ is the vector field defined on $\Sigma$, as the
first-order variation of the unit normal vector field.
\end{lemma}

\begin{proof}
Set $Y:=\overline{\nabla}_XN$.
We need the  following well known formula, of which we include a proof for the 
reader's convenience:
\be \label{eq:div}
 div(T_p Y) & = & S'_{p+1} + \frac{1}{2}\langle I', BT_p\rangle +
(n-p)K S_p\langle X, N\rangle~.
\ee
Indeed, consider a deformation $\phi_t,\ t\in(-\epsilon,\epsilon),\
\epsilon>0$,  of the immersion 
$\phi$ having $X$ as deformation vector field. Denote by $N_t$ a unit normal
field to $\phi_t$ 
depending in a differentiable way on $t$. Call $I_t$  the  metric
induced on $\Sigma$ 
by $\phi_t$ and $B_t$ its shape operator.
Let $u,v$ be tangent vectors to $\Sigma$. Taking the derivative at $t=0$ in
the equation: 
$$I_t(B_tu,v) = \left\langle \overline{\nabla}_{d\phi_t(u)} N_t ,
d\phi_t(v)\right\rangle,$$ 
we obtain:
$$ I'(Bu,v)+ I(B'u,v)=
\left\langle \left(\overline{\nabla}_{\frac{\partial\phi_t}{\partial
    t}}\overline{\nabla}_{d\phi_t(u)} N_t\right)_{|t=0} , v \right\rangle 
+ \langle  \overline{\nabla}_u N, \overline{\nabla}_v 
X\rangle. $$
Now, on the one hand we have
$$ I'(Bu,v)= \langle \overline{\nabla}_{Bu}X, v\rangle + \langle Bu,
\overline{\nabla}_v X\rangle~. $$ 
On the other hand, denoting by $R$ the curvature tensor of $M_K^{n+1}$, we can
write: 
$$ \left\langle \left(\overline{\nabla}_{\frac{\partial\phi_t}{\partial
    t}}\overline{\nabla}_{d\phi_t(u)} N_t\right)_{t=0} , v\right\rangle 
=\langle \overline{\nabla}_{u}\overline{\nabla}_{X} N , v\rangle - \langle
R(X,u)N,v\rangle~. $$  
Collecting the terms we finally get:
\be\label{eq:der}
 I(B'u,v)= \langle \overline{\nabla}_u Y, v\rangle - \langle
\overline{\nabla}_{Bu} X, v\rangle  - \langle R(X,u)N,v\rangle~. 
\ee
Take now a local orthonormal frame $\{e_1,\ldots,e_n\}$ on $\Sigma.$ Then:
\begin{eqnarray*}
div(T_p Y) & = & \sum_{i=1}^n \langle \nabla_{e_i}(T_pY), e_i\rangle \\
&=&  \sum_{i=1}^n \langle (\nabla_{e_i}T_p)Y +
T_p \nabla_{e_i}Y , e_i\rangle \\
& = & div(T_p)( Y)
+ \sum_{i=1}^n \langle \nabla_{e_i}Y , T_p e_i\rangle \\
&=& \sum_{i=1}^n \langle \nabla_{e_i}Y , T_p e_i\rangle .
\end{eqnarray*}
We now use equation (\ref{eq:der}) above to write:
\begin{eqnarray*}
div(T_p Y) & = & \sum_{i=1}^n  I(B^{\prime} e_i, T_p e_i) +  \langle
\overline{\nabla}_{Be_i} X, e_i\rangle  + \langle R(X,e_i)N,e_i\rangle \\
&=& tr (B^\prime T_p) + \frac{1}{2}\langle I', BT_p\rangle +K \sum_{i=1}^n\langle X,N\rangle\langle e_i,T_pe_i\rangle\\
&=& tr (B^\prime T_p) + \frac{1}{2}\langle I', BT_p\rangle +K \langle X,N\rangle tr(T_p) .
\end{eqnarray*}
Finally using properties (2) and (4) of Lemma \ref{newton} we obtain the desired equation (\ref{eq:div}) above. 

Now, for any fixed vector $a\in \R^{n+2}:$
$$div(\langle a,x\rangle T_pY) = \langle \nabla(\langle x,a\rangle),
T_pY \rangle + \langle a,x\rangle div(T_pY),$$
and so using (\ref{eq:div}) we obtain:
\be \label {eq:divY}
div(\langle a,x\rangle T_pY)= \langle a,T_pY\rangle + \langle a,x\rangle \left( S'_{p+1} +
\frac{1}{2}\langle I', BT_p\rangle + (n-p)K S_p\langle X, N\rangle\right)~,
\ee
and the result follows by integration over $\Sigma$.
\end{proof}
We will see (\ref{sphypset}) that this argument can be applied when $\Sigma$ is an (outer) $\epsilon-$neighborhood of a polyhedron, although it is only $C^{1,1}$ rather than $C^2.$

\section{Polyhedral localization formulas}\label{sec:localization}

This section contains some technical statements on the limits of integral
quantities defined on a hypersurface, when one considers them on the boundary
of the  (outer) $\epsilon$-neighborhood of a polyhedron $P$, as $\epsilon\rightarrow
0$. The general idea is that some integral quantities have limits that depend
only on what happens  on faces of $P$ of given codimension.
For any $\epsilon>0$ small enough, we call $P_\epsilon$ the boundary of the
set of points at distance at most $\epsilon$ from the interior of $P$.
Given a face $F$ of $P$, of codimension $p+1$, we call $F_\epsilon$ the set of
points $x\in P_\epsilon$ such that the point of $P$ which is closest to $x$ is
in $F$.  We orient $P$ and $P_{\epsilon}$ by the exterior unit normal.

\subsection{The spherical case.}

We now concentrate on the case where the ambient space is $S^{n+1}$. It turns
out that the same results apply in the hyperbolic space, and the proofs
can be extended from the spherical to the hyperbolic context with only
very little changes.
The local geometry of $F_\epsilon$ is simple. Its shape operator has
two eigenvalues: 
\begin{itemize}
\item $1/\tan(\epsilon)$, which has multiplicity $p$, on directions
corresponding (by parallel transport along shortest geodesics) to directions
orthogonal to $F$.
\item $- \tan(\epsilon)$, which has multiplicity $n-p$, on directions
corresponding by parallel transport to directions tangent to $F$.
\end{itemize}
Actually $F_\epsilon$ is isometric to the product $\cos(\epsilon)F\times
\sin(\epsilon) F^*$, so that its volume is $\cos(\epsilon)^{n-p}
\sin(\epsilon)^p V(F) V(F^*)$.

\begin{remark}
$$ \lim_{\epsilon\rightarrow 0}\int_{P_\epsilon} S'_p x dv = 0~. $$
\end{remark}

\begin{proof}
This is clear since, for each face $F$ of $P$, $S'_p=0$ on $F_\epsilon$.
\end{proof}

We now consider in more details how some important quantities behave on the
surfaces $F_\epsilon$ which ``approximates'' a codimension $p+1$ face of $P$.
At each point of such a surface, the tangent space $T_xF_\epsilon$ has a natural
decomposition as the direct sum of two subspaces, corresponding under the
parallel transport along the geodesic segments orthogonal to $F_\epsilon$ to
the spaces tangent (resp. orthogonal) to $F$. We call those subspaces
$V$ (resp. $W$) and also call $\Pi_V$ (resp. $\Pi_W$) the orthogonal
projection from $T_xF_\epsilon$ to $V$ (resp. to $W$). The dimension of
$W$ is equal to $p$, while the dimension of $V$ is equal to the dimension
of $F$, which is $n-p$ since $F$ has codimension $p+1$.

\begin{lemma} \label{lm:Tp}
Let $F\in F_{p+1}$ be a codimension $p+1$ face.
For all $q\in \{ 1,\cdots, n\}$, $S_q$ is constant over $F_\epsilon$.
Moreover:
\begin{enumerate}
\item for  $q\neq p$, $S_q=o(\epsilon^{-p})$ as
$\epsilon\rightarrow 0$, while $S_p=\epsilon^{-p}+o(\epsilon^{-p})$,
\item $T_q = o(\epsilon^{-p})$ if $q\neq p$, while $T_p = \epsilon^{-p}
\Pi_V + o(\epsilon^{-p})$,
\item $BT_q = o(\epsilon^{-p})$ for $q\neq p-1$, while
$BT_{p-1} = \epsilon^{-p} \Pi_W + o(\epsilon^{-p})$.
\end{enumerate}
\end{lemma}

\begin{proof}
At each point of $F_\epsilon$ the shape operator has only two eigenvalues
(corresponding to two different principal curvatures): $-\tan(\epsilon)$,
with eigenspace $V$, and $1/\tan(\epsilon)$ with eigenspace $W$.
The result concerning $S_p$ then follows from the explicit expression
of $S_p$, since at most $p$ principal curvatures can be equal to
$1/\tan(\epsilon)$ while the others are equal to $-\tan(\epsilon)$.
To obtain the result on $T_p$ and on
$BT_p$ it is necessary to remark that, for each
eigenvector $e_i$ of the shape operator at a point $x\in F_\epsilon$, $e_i$ is
also an eigenvector of $T_p$ with eigenvalue equal to the $p$-th symmetric
function in the other eigenvalues of $B$ at $x$:
$$ T_p e_i = \left( \sum_{\stackrel{i_1<i_2<\cdots <i_p}{i_k\neq i,
1\leq k\leq p}}
k_{i_1}k_{i_2}\cdots k_{i_p} \right) e_i~. $$
This formula, which can be checked by a direct induction argument using the
  definition of the Newton operators, leads directly to the second result:
$T_q$ is maximal for $q=p$, and is then equivalent to $\epsilon^{-p} \Pi_V$.
Similarly, $BT_q$ is
  maximal for $q=p-1$, and is then equivalent to $\epsilon^{-p}\Pi_W$.
\end{proof}

It is now possible to estimate how some interesting quantities behave as
$\epsilon\rightarrow 0$. We introduce for this a simple notation. 

\begin{df}
Consider a first-order deformation of the polyhedron $P$ with deformation
field a vector field 
$X$ tangent to $S^{n+1}$ defined on $P$. Let $F$ be a face of $P$. For each
point $x\in F$, we call $\nu_F(x)$ the 
component of $X$ at $x$ which is orthogonal to the subspace of $\R^{n+2}$
generated by $F$, called $vect(F)$ here. Then $\nu_F$ is the restriction
to $F$ of a linear map, which we still call $\nu_F$, from $vect(F)$ to
$vect(F^*)$. 
\end{df}

\begin{remark}
The adjoint operator $\nu_F^*:vect(F^*)\rightarrow vect(F)$ is equal to
$-\nu_{F^*}$.
\end{remark}

\begin{proof}
Let $x\in F, y\in F^*$, then $\langle x,y\rangle=0$, and this remains true
under a first-order deformation. So, under this deformation
$$
\langle x,y\rangle' = \langle \nu_F(x),y\rangle + \langle x,
\nu_{F^*}(y)\rangle = 0~, 
$$
and the result follows.
\end{proof}

\begin{lemma} \label{lm:limS}
For all $p\in \{ 0,\cdots, n-1\}$:
$$ \lim_{\epsilon\rightarrow 0}\int_{P_\epsilon} fS_p x dv = \sum_{F\in
F_{p+1}} \int_F \langle \nu_F(x), \pi(F^*)\rangle x dv~. $$
\end{lemma}

\begin{proof}
Clearly $f$ is bounded over $P_\epsilon$, and the bound is uniform as
$\epsilon\rightarrow 0$. We have seen above that, given a face $F$
of $P$ of codimension $q+1$, the volume of $F_\epsilon$ is equal to
$\sin(\epsilon)^{q} \cos(\epsilon)^{n-q}V(F)V(F^*)$.   
Given the estimates on $S_p$ described above, it follows that:
$$ \int_{F_\epsilon} fS_p x dv \rightarrow 0 $$
unless $q=p$, and that in that case, since $f:=\langle X,N\rangle$:
$$ \int_{F_\epsilon} fS_p x dv \rightarrow \int_{x\in F, n\in F^*}
\langle n,X\rangle x dv(x) dv^*(n)~, $$
where $dv^*$ is used to denote the volume element on $F^*$.
But by definition of $\nu_F$ we have: $\langle n,X\rangle =
\langle\nu_F(x),n\rangle$,
at each point $x\in F$ and for each normal vector $n\in F^*$, and
it follows that:
$$ \int_{F_\epsilon} fS_p x dv \rightarrow \int_{x\in F, n\in F^*}
\langle \nu_F(x),n\rangle x dv(x) dv^*(n) = \int_{x\in F} \langle \nu_F(x),
\pi(F^*)\rangle xdv(x)~, $$
and the result follows by summing over faces of codimension $p+1$.
\end{proof}
It will be useful below to note that this formula extends to deformations
of polyhedra which, rather than being contained in the ``usual'' unit sphere
$S^{n+1}\subset \R^{n+2}$, are contained in a unit sphere centered at another
point $a\in \R^{n+2}$.
In that case, call $\tau_{-a}$ the translation by $-a$ in $\R^{n+2}$,
and let $P_0:=\tau_{-a}(P)$, so that $P_0\subset S^{n+1}$.
For each face $F$ of $P$, define
the dual $F^*$ of $F$ as $\tau_{-a}(F)^*$ (which is
well-defined since $\tau_{-a}(F)$ is a face of $P_0\subset S^{n+1}$).
It is also necessary to replace the function $\nu_F$ defined
above (for any face $F$ of $P$) by another function, which we call
$\overline{\nu}_F$, and which is defined as $\nu_{\tau_{-a}(F)}\circ\tau_{-a}$,
where $\nu_{\tau_{-a}(F)}$ is defined with respect to the first-order
deformation of $P_0$ corresponding to the chosen deformation of $P$.

\begin{remark} \label{rk:34}
With those notations, the formula becomes:
$$ \lim_{\epsilon\rightarrow 0}\int_{P_\epsilon} fS_p x dv = \sum_{F\in
F_{p+1}} \int_F \langle \overline{\nu}_F(x), \pi(F^*)\rangle x dv~. $$
\end{remark}

\begin{proof}
The proof can be obtained as in the previous lemma.
Alternatively it is possible to apply the previous lemma to $\tau_{-a}P$,
using the fact that:
$$ \int_{P_\epsilon} fS_p xdv = \int_{\tau_{-a}(P_\epsilon)} fS_pxdv +
a\int_{\tau_{-a}(P_\epsilon)} fS_p dv~, $$
while:
$$ \sum_{F\in F_{p+1}}
\int_F \langle \overline{\nu}_F(x), \pi(F^*)\rangle x dv
= \sum_{F\in F_{p+1}}
\int_{\tau_{-a}(F)}
\langle \overline{\nu}_F(x), \pi(F^*)\rangle x dv + a\sum_{F\in F_{p+1}}
\int_{\tau_{-a}(F)}\langle \overline{\nu}_F(x), \pi(F^*)\rangle dv~, $$
The result then follows from the fact that, for any polyhedron in $S^{n+1}$:
$$ \lim_{\epsilon\rightarrow 0}\int_{P_\epsilon} fS_p dv = \sum_{F\in
F_{p+1}} \int_F \langle \overline{\nu}_F(x), \pi(F^*)\rangle dv~, $$
which can be proved exactly like the previous lemma (removing some
``$x's$'' from the equations).
\end{proof}

\begin{lemma} \label{lm:35}
For any first-order deformation of $P$ and any $p\in \{ 0, \cdots, n-1\}$:
$$ \lim_{\epsilon\rightarrow 0}\int_{P_\epsilon} \frac{1}{2}\langle I',
BT_p\rangle x dv = \sum_{F\in F_{p+2}} V'(F^*) \pi(F)~. $$
\end{lemma}

\begin{proof}
The same argument as in the proof of Lemma \ref{lm:limS} shows that the
contribution to the limit of faces of codimension $q$ vanishes unless
$q=p+2$. For faces of codimension $p+2$, $BT_p\simeq \epsilon^{-p-1} \Pi_W$,
so that:
$$ \frac{1}{2}\langle I', BT_p\rangle dv^* \simeq \epsilon^{-p-1} (dv^*)'~, $$
where $(dv^*)'$ is the first-order variation of the volume element on $V^*$.
Therefore:
$$ \lim_{\epsilon\rightarrow 0}\int_{P_\epsilon} \frac{1}{2}\langle I',
BT_p\rangle x dv = \sum_{F\in F_{p+2}} \int_{x\in F, n\in F^*}
x dv(x) (dv^*)'(n) = \sum_{F\in F_{p+2}} V'(F^*) \pi(F)~. $$
\end{proof}

\begin{remark}
The same formula is valid if one considers a sphere centered at a
point $a\in \R^{n+2}$ which is not the origin.
\end{remark}

\begin{proof}
Again the proof can be obtained either by checking that the proof of the
previous lemma extends to a sphere centered at a point $a\in \R^{n+2}$, or
by remarking that the difference in the two terms that arises by centering
the sphere at $a$ is equal to:
$$ a\left(\lim_{\epsilon\rightarrow 0}\int_{P_\epsilon}\frac{1}2\langle I',
BT_p\rangle dv -  \sum_{F\in F_{p+2}} V'(F^*) V(F)\right)~, $$
a term which is shown to vanish by the arguments used to prove the previous
lemma, with some ``$x$'' suppressed from the equations.
\end{proof}

\begin{lemma} \label{lm:37}
Again for any first-order deformation of $P$ and for any $p\in \{0,\cdots,
n-1\}$: 
$$ \lim_{\epsilon\rightarrow 0}\int_{P_\epsilon} T_p 
\overline{\nabla}_XN dv =
\sum_{F\in F_{p+1}} \left(- V(F)\nu_F^*(\pi(F^*))
+ K \int_F \langle \nu_F(x), \pi(F^*)\rangle x dv\right)~. $$
\end{lemma}

\begin{proof}
First note that $\overline{\nabla}_XN$ is bounded, so that,
by the second point of Lemma \ref{lm:Tp}, only faces
of codimension $p+1$ give a non-zero contribution to the limit on the
left-hand side. It follows that:
$$ \lim_{\epsilon\rightarrow 0}\int_{P_\epsilon} T_p \overline{\nabla}_XN dv =\lim_{\epsilon\rightarrow 0}
\sum_{F\in F_{p+1}} \int_{F_\epsilon}
\Pi_{V}(\overline{\nabla}_XN)dv~. $$
This can also be written as:
\be \label{eq:diff1}
\lim_{\epsilon\rightarrow 0}\int_{P_\epsilon} T_p \overline{\nabla}_XN dv & = &
\sum_{F\in F_{p+1}}\int_F\int_{F^*} \Pi_{T_xF} (\nu_{F^*}(n)) dv(x) dv^*(n)~.
\ee
Note however that $\nu_{F^*}(n)\in vect(F)$, so that 
$$ \Pi_{T_xF} (\nu_{F^*}(n)) = \nu_{F^*}(n) - K\langle \nu_{F^*}(n), x\rangle
x~, $$
where the term ``$K$'' is necessary because, if $K=-1$, $x$ is a time-like
unit vector in the Minkowski space $\R^{n+2}_1$.
So it follows from (\ref{eq:diff1}) that:
\be \label{eq:diff}
\lim_{\epsilon\rightarrow 0}\int_{P_\epsilon} T_p \overline{\nabla}_XN dv & = &
\sum_{F\in F_{p+1}} \int_F\int_{F^*} \nu_{F^*}(n) - K \langle
\nu_{F^*}(n),x\rangle x dv(x) dv^*(n)~.
\ee
After integration in $n$ this becomes:
$$ \lim_{\epsilon\rightarrow 0}\int_{P_\epsilon} T_p \overline{\nabla}_XN dv =
V(F) \nu_{F^*}(\pi(F^*)) - K \int_F \langle \nu_{F^*}(\pi(F^*)),x\rangle xdv~. $$
We have seen above that $\nu_{F^*}=-(\nu_F)^*$. So:
$$ \lim_{\epsilon\rightarrow 0}\int_{P_\epsilon} T_p \overline{\nabla}_XN dv =
- V(F) \nu_F^*(\pi(F^*)) + K \int_F \langle \pi(F^*), \nu_F(x)\rangle xdv~,$$
and the result follows.
\end{proof}

\begin{remark}
Consider a unit sphere centered at a point $a\in \R^{n+2}$. The formula
becomes:
$$ \lim_{\epsilon\rightarrow 0}\int_{P_\epsilon} T_p \overline{\nabla}_XN dv =
\sum_{F\in F_{p+1}} \left(- V(F)\nu_F^*(\pi(F^*))
+ \int_F \langle \overline{\nu}_F(x), \pi(F^*)\rangle x dv -
\langle\int_F\nub_F(x)dv, \pi(F^*)\rangle a\right)~. $$
\end{remark}

\begin{proof}
The proof proceeds as above, except that equation (\ref{eq:diff})
is replaced by:
$$ \lim_{\epsilon\rightarrow 0}\int_{P_\epsilon} T_p\overline{ \nabla}_XN dv =
\sum_{F\in F_{p+1}} \int_F\int_{F^*} \nu_{F^*}(n) - \langle
\nu_{F^*}(n),(x-a)\rangle (x-a) dv^*(n) dv(x)~. $$
Therefore:
$$ \lim_{\epsilon\rightarrow 0}\int_{P_\epsilon} T_p\overline{\nabla}_XN dv =
V(F) \nu_{F^*}(\pi(F^*)) -
\int_F \langle \nu_{F^*}(\pi(F^*)),x-a\rangle (x-a)dv~. $$
So:
$$ \lim_{\epsilon\rightarrow 0}\int_{P_\epsilon} T_p\overline{\nabla}_XN dv =
- V(F) \nu_F^*(\pi(F^*)) +
\int_F \langle \pi(F^*), \overline{\nu}_F(x)\rangle (x-a)dv~,$$
and the result follows.
\end{proof}

\begin{lemma} \label{lm:39}
For all $p\in \{ 0,\cdots, n-1\}$:
$$ \lim_{\epsilon\rightarrow 0}\int_{P_\epsilon} T_p \xi + \frac{x}{2} \langle
I', T_p\rangle dv = \sum_{F\in
F_{p+1}} V(F^*)(\pi'(F) - \nu_F(\pi(F)))~. $$
\end{lemma}

\begin{proof}
The same argument as above can be used to show that only codimension $p+1$
faces have a non-zero contribution to the limit. So:
\begin{eqnarray*}
\lim_{\epsilon\rightarrow 0}\int_{P_\epsilon} T_p \xi + \frac{x}{2} \langle
I', T_p\rangle dv & = &
\sum_{F\in F_{p+1}} V(F^*)\int_F\Pi_{T_xF}(\xi) + \frac{x}{2} \langle I',
I\rangle dv \\
& = &
\sum_{F\in F_{p+1}} V(F^*)\int_F\Pi_{T_xF}(X) + \frac{x}{2} \langle
I', I\rangle dv \\
& = &
\sum_{F\in F_{p+1}} V(F^*)\int_F X - \nu_F(x) + \frac{x}{2} \langle
I', I\rangle dv \\
& = & \sum_{F\in F_{p+1}} V(F^*) \left( \left(\int_F xdv\right)' -
\int_F \nu_F(x)dx \right) \\
& = & \sum_{F\in F_{p+1}} V(F^*) (\pi'(F) -
\nu_F(\pi(F)))~,
\end{eqnarray*}
which is the desired result.
\end{proof}

\begin{remark} \label{rk:310}
Consider a sphere centered at a point $a$ which is not the origin. Then  the formula
becomes:
$$ \lim_{\epsilon\rightarrow 0}\int_{P_\epsilon} T_p \xi + \frac{x}{2} \langle
I', T_p\rangle dv = \sum_{F\in
F_{p+1}} V(F^*)\left(\pi'(F) - \int_F\nub_F(x)\right)dv~. $$
\end{remark}

\begin{proof}
The proof follows directly from the computations above.
\end{proof}

\section{Proof of the main results}

In this section we use the results of section 4 to take the limit of 
the smooth formulas of section 3 for (outer) $\epsilon$-neighborhoods of 
polyhedra, to obtain the proofs of the main results. 

\subsection{The spherical and hyperbolic settings.} \label{sphypset}

We first explain why  the formulas in Lemma \ref{lm:23} and Lemma \ref{lm:24} do apply to outer $\epsilon-$neighborhoods ($\epsilon>0$ small enough) 
of convex polyhedra although they are only $C^{1,1}.$ Indeed,  $P_{\epsilon}$ can be decomposed  as 
$$P_{\epsilon}=\bigcup_{{m=1} \atop{F\in F_m}}^{m=n+1} F_{\epsilon},$$
where $F_{\epsilon}$ denotes, as before,  the set of points where the normal to $P_{\epsilon}$ meets $P$
at the  codimension $m$ face $F.$ The parts $F_{\epsilon}$ are smooth and have piecewise smooth boundaries. We can thus apply to each $F_{\epsilon}$ the formulas (\ref{eq:divxi}) and (\ref{eq:divY}) in the proofs of Lemmas
\ref{lm:23} and \ref{lm:24}, respectively. 
In order to show that the formulas of Lemma \ref{lm:23} and Lemma \ref{lm:24}
apply to $P_{\epsilon}$, we have to check that the boundary terms that
appear after integrating, and using the divergence theorem,  in (\ref{eq:divxi}) and (\ref{eq:divY}), cancel out two by two after
summing up. These boundary terms are
(using the symmetry of the Newton transformations) of the form
$\int_{\partial F_{\epsilon}}\langle a,x\rangle \langle Z, T_r(\nu)\rangle$, for some
$1\le r\le n$, $Z$ being
a (continuous) vector field on $P_{\epsilon}$ and $\nu$ the unit exterior
conormal to
$\partial F_{\epsilon}$. Let  $F\in F_m$ and $Q\in F_q,$ then $F_{\epsilon}$ and $ Q_{\epsilon}$
have a common boundary of non-zero $(n-1)$-measure if and only if $q=m+1$
or (symmetrically) $q=m-1$. Take then $F\in F_m$ and $Q\in F_{m+1}.$ Along the common boundary of $F_{\epsilon}$
and $Q_{\epsilon}$, the  unit conormal $\nu$ corresponds under parallel transport along geodesic segments orthogonal to $F_{\epsilon}$ (resp. $Q_{\epsilon}$) to the   space tangent to $F$ (resp. 
orthogonal to $Q$). Then a straightforward computation, using properties (1), (2) and (3) in Lemma \ref{lm:Tp},
shows  that the vectors $T_r(\nu)$, where $\nu$ is a vector orthogonal to
the common
boundary, are equal on
the two sides of the common boundary (cf. \cite{hsf}).

We can now state two lemmas obtained by applying the polyhedral localization
formulas to Lemma \ref{lm:23} and Lemma \ref{lm:24}. 
\begin{lemma} \label{lm:41}
Let $P\subset M_K^{n+1}, K\in \{-1,1\}$, and let $X$ be a first-order
deformation vector field on $P$. For all $p\in \{ 1,2,\cdots, n\}$:
\begin{equation*} \tag{\mbox{$F'_p$}} 
\sum_{F\in F_p} V(F^*) (\pi'(F) - \nu_F(\pi(F))) - p \sum_{G\in F_{p+1}} \int_G
\langle \nu_G(x), \pi(G^*)\rangle x dv = 0~. 
\end{equation*}
\end{lemma}

\begin{proof}
According to Lemma \ref{lm:23}, we have for all $p\in \{ 0,\cdots, n-1\}$:
\begin{equation*}
\int_{P_{\epsilon}} T_p \xi + \frac{x}{2} \langle I', T_p\rangle -  (p+1) xfS_{p+1}
dv = 0~.
\end{equation*}
However we know by Lemma
\ref{lm:limS} that 
$$ \lim_{\epsilon\rightarrow 0}\int_{P_\epsilon} fS_{p+1} x dv = \sum_{F\in
F_{p+2}} \int_F \langle \nu_F(x), \pi(F^*)\rangle x dv~, $$
while Lemma \ref{lm:39} shows that
$$ \lim_{\epsilon\rightarrow 0}\int_{P_\epsilon} T_p \xi + \frac{x}{2} \langle
I', T_p\rangle dv = \sum_{F\in
F_{p+1}} V(F^*)(\pi'(F) - \nu_F(\pi(F)))~. $$
The proof follows by replacing $p$ by $p-1$. 
\end{proof}

\begin{remark} \label{local:a1} 
If the sphere is centered at a point $a\in
  \R^{n+2}$ then the formula becomes: 
$$ \sum_{F\in F_p} V(F^*) (\pi'(F) - \int_F\nub_{F}(x)dv) - 
p \sum_{G\in F_{p+1}} \int_G 
\left\langle \overline{\nu}_G(x), \pi(G^*)\right\rangle x dv = 0~. $$
\end{remark}

\begin{proof}
As the proof of the previous Lemma, with Remark \ref{rk:34} instead of
Lemma \ref{lm:limS} and Remark \ref{rk:310} rather than Lemma \ref{lm:39}. 
One should then remark that, by the definition of $\nu_H$,
$$ \nu_H^*(\pi(H^*)) = -\Pi_H(\pi'(H^*))~. $$
\end{proof}

\begin{lemma} \label{lm:43}
For all $p=0,\cdots, n-1$:
\begin{equation*} \tag{\mbox{$G'_p$}}
\sum_{F\in F_{p+2}} V'(F^*) \pi(F) + (n - p + 1)K \sum_{G\in F_{p+1}} \int_G
\langle \nu_G(x), \pi(G^*)\rangle x dv - \sum_{H\in F_{p+1}} V(H)
\nu_H^*(\pi(H^*))  = 0~.
\end{equation*}
\end{lemma}

\begin{proof}
Recall that by Lemma \ref{lm:24} above,
\begin{equation*}
\int_{P_\epsilon} \left(
S'_{p+1} + \frac{1}{2} \langle I', BT_p\rangle + (n-p) Kf S_p\right)
x + T_p  \overline{\nabla}_XN dv = 0~.
\end{equation*}
However, by Lemma \ref{lm:limS},
$$ \lim_{\epsilon\rightarrow 0}\int_{P_\epsilon} fS_p x dv = \sum_{F\in
F_{p+1}} \int_F \langle \nu_F(x), \pi(F^*)\rangle x dv~, $$
while Lemma \ref{lm:35} shows that
$$ \lim_{\epsilon\rightarrow 0}\int_{P_\epsilon} \frac{1}{2}\langle I',
BT_p\rangle x dv = \sum_{F\in F_{p+2}} V'(F^*) \pi(F)~. $$
and Lemma \ref{lm:37} indicates that
$$ \lim_{\epsilon\rightarrow 0}\int_{P_\epsilon} T_p \overline{\nabla}_XN dv =
\sum_{F\in F_{p+1}} \left(- V(F)\nu_F^*(\pi(F^*))
+ K \int_F \langle \nu_F(x), \pi(F^*)\rangle x dv\right)~. $$
The proof follows.
\end{proof}

\begin{remark}\label{local:a2} 
If the sphere is centered at a point $a\in
  \R^{n+2}$ then the formula becomes: 
$$  \sum_{F\in F_{p+2}} V'(F^*) \pi(F) + (n-p+1) K \sum_{G\in F_{p+1}} \int_G
\langle \overline{\nu}_G(x), \pi(G^*)\rangle x dv - 
\coupeq - \sum_{H\in F_{p+1}}
\left(V(H) 
\nu_H^*(\pi(H^*))  +\langle\int_H\nub_H(x)dv, \pi(H^*)\rangle\
a \right) = 0~. $$ 
\end{remark}

\begin{proof}[Proof of Theorem \ref{tm:general}]
Formula ($E_p$) is obtained as the sum of $(n-p+1)K$ times formula ($F'_p$)
in Lemma \ref{lm:41} and $p$ times formula ($G'_p$) in Lemma \ref{lm:43}.
\end{proof}

\begin{proof}[Proof of Theorem \ref{tm:codim1}]
Lemma \ref{lm:43}, applied with $p=0$, translates as:
$$  \sum_{F\in F_{2}} \theta'_F \pi(F) + (n+1) K \sum_{G\in F_{1}} \int_G
\langle \nu_G(x), N_G\rangle x dv +\sum_{H\in F_{1}} V(H)
N'_H  = 0~. $$
But it is easy to check that:
$$ \left(\int_{int(P)} xdv\right)' = \sum_{G\in F_1} \int_G
\langle \nu_G(x), N_G\rangle x dv~, $$
so that:
$$  \sum_{F\in F_{2}} \theta'_F \pi(F) + (n+1) K \left(\int_{int(P)} xdv\right)'
+ \sum_{H\in F_{1}} V(H)N'_H = 0~. $$
Using the first equation in Remark \ref{rk:18} leads to:
$$  \sum_{F\in F_{2}} \theta'_F \pi(F) - \left(\sum_{H\in F_{1}} V(H)N_H\right)'
+ \sum_{H\in F_{1}} V(H)N'_H = 0~, $$
and the result follows.
\end{proof}

\subsection{The second main formula.}

The proof of Theorem \ref{tm:quatre} follows from establishing a formula which is {\it dual} to the one in Theorem
\ref{tm:general} (see the Remark below). This is the content of the following

\begin{prop}\label{dual} In any first-order deformation of a convex polyhedron $P$ in $M^{n+1}_K$ and for 
any $p\in \{ 1, 2, \cdots, n-1\}$:
\begin{equation*}\tag{\mbox{$L_p$}} 
(p+1) \sum_{F\in F_{p+2}} V(F)(\pi'(F^*)-\nu_{F^*}(\pi(F^*))+ K(n-p)\sum_{G\in F_{p}} V^\prime(G) \pi(G^*)=
-K(n-p)\sum_{H\in F_{p+1}} V(H^*)\nu_H(\pi(H))
~. \end{equation*}
\end{prop}

\begin{proof} The idea is to replace  in the proof of Theorem \ref{tm:general} the position vector $x$ by the (exterior) normal vector $N.$  We will indicate only the main steps in the proof and omit the details which are very similar 
to the arguments involved in the proof of Theorem \ref{tm:general}.

Consider a compact hypersurface $\Sigma$ in $M^{n+1}_K,$ oriented by a unit normal field $N.$ Consider a first order deformation of $\Sigma$ with deformation vector field $X$ and let $X= \xi+ fN$ be its decomposition into its tangent part  $\xi$
to  $\Sigma$  and normal one $fN.$ For any fixed vector $a,$ in $\R^{n+2}$ in case $K=1$ and in $\R^{n+2}_1$ in case $K=-1,$ and any 
$p\in \{o,\ldots,n-1\},$ we have: 
$$ div(\langle a,N\rangle T_p\xi) = \left\langle a, \left(  \frac{1}{2}\langle I',T_p\rangle -(p+1)fS_{p+1}\right)N + BT_p\xi\right\rangle.$$
Integrating over $\Sigma$ leads to the equation, valid for $p=0, \cdots, n-1$:
\begin{equation*}\tag{\mbox{$Q_p$}}
\int_\Sigma BT_p \xi + \frac{N}{2} \langle I', T_p\rangle -  (p+1) NfS_{p+1}
dv = 0~.
\end{equation*}

The localization technique as in Section \ref{sec:localization} gives the following formula, valid for  all $p=0,\cdots, n-1,$ for first order deformations of a convex polyhedron $P\subset M^{n+1}_K:$
\begin{equation*} \tag{\mbox{$Q'_p$}}
\sum_{F\in F_{p+2}} V(F^*)\nu_F( \pi(F) ) - (p+2) \sum_{F\in F_{p+2}} \int_{F^*}
\langle \nu_F(\pi(F)), n\rangle n dv^*(n)  + \sum_{G\in F_{p+1}} V^\prime(G)\pi(G^*) = 0~.
\end{equation*}

Next,  we have for all $p=0,\cdots, n-1:$
$$
div(\langle a,N\rangle T_p\overline{\nabla}_XN)= \left\langle a, \left( S'_{p+1} +  \frac{1}{2}\langle I',BT_p\rangle 
+(n-p)KfS_p\right)N + BT_p\xi\right\rangle.$$
Integrating over $\Sigma$ leads to the equation
\begin{equation*}\tag{\mbox{$M_p$}}
\int_\Sigma \left(
S'_{p+1} + \frac{1}{2} \langle I', BT_p\rangle + (n-p) Kf S_p\right)
N+ BT_p \overline{\nabla}_XN dv = 0~,
\end{equation*}
Again, the localizaion technique gives that,
 for all $p\in \{ 0,\cdots, n-1\}$:
\begin{equation*} \tag{\mbox{$M'_p$}} 
\sum_{F\in F_{p+2}} V(F) (\pi'(F^*) - \nu_{F^*}(\pi(F^*))) +(n-p)K  \sum_{G\in F_{p+1}} \int_{G^*}
\langle \nu_G(\pi(G)), n\rangle n dv^*(n) = 0~. 
\end{equation*}
Formula $(L_p)$ follows taking, for $p\in\{1,\ldots,n-1\},$  the combination $K(n-p) (Q'_{p-1})+(p+1)(M'_{p}).$
\end{proof}

\begin{proof}[Proof of Theorem \ref{tm:quatre}]
The proof is obtained adding, for $p\in \{  2, \cdots, n-1\},$ 
equation $(E_p)$ in Theorem \ref{tm:general} to equation $(L_{p-1})$ 
in Proposition \ref{dual} above and using 
the fact that $\nu^*_G=-\nu_{G^*}.$ 
\end{proof}
 
\begin{remark}
In the spherical case,  formula $(L_p)$ in Proposition \ref{dual} 
follows in fact directly under duality from Theorem 
\ref{tm:general}. 
\end{remark}

\begin{proof}
 Indeed, let $p\in \{ 1,\cdots, n-1\}$. By Theorem \ref{tm:general}:
$$
(n-p+1)\sum_{H\in F_p} V(H^*)(\pi'(H) - \nu_H(\pi(H))) 
+ p\sum_{F\in F_{p+2}} V'(F^*) \pi(F) = p \sum_{G\in F_{p+1}} V(G)
\nu_G^*(\pi(G^*))~. 
$$
We can apply this formula to the dual polyhedron $P^*$ for an integer 
$q\in \{ 1, \cdots, n-1\}$. We call $F^*_l$ the set of faces of codimension
$l$ of $P^*$, and find that: 
$$ (n-q+1)\sum_{H^*\in F^*_q} V(H)(\pi'(H^*) - \nu_{H^*}(\pi(H^*))) 
+ q\sum_{F^*\in F^*_{q+2}} V'(F) \pi(F^*) = q \sum_{G^*\in F^*_{q+1}} V(G^*)
\nu^*_{G^*}(\pi(G))~. $$
But a face of codimension $l$ of $P^*$ is the dual of a face of codimension
$n+2-l$ of $P$. Moreover, we have seen that $\nu^*_{G^*}=-\nu_G$. So the
formula above can be written as follows:
$$ (n-q+1)\sum_{H\in F_{n+2-q}} V(H)(\pi'(H^*) - \nu_{H^*}(\pi(H^*))) 
+ q\sum_{F\in F_{n-q}} V'(F) \pi(F^*) = - q \sum_{G\in F_{n+1-q}} V(G^*)
\nu_{G}(\pi(G))~. $$
Take $p\in \{ 2, \cdots, n\}$, then $n+1-p\in \{ 1, \cdots, n-1\}$, so that it
is possible to set $q=n+1-p$, so that the equation becomes:
$$ p\sum_{H\in F_{p+1}} V(H)(\pi'(H^*) - \nu_{H^*}(\pi(H^*))) 
+ (n+1-p)\sum_{F\in F_{p-1}} V'(F) \pi(F^*) = 
- (n+1-p) \sum_{G\in F_{p}} V(G^*)
\nu_{G}(\pi(G))~. $$
Replacing $p$ by $p+1$ leads to $(L_p).$
\end{proof}

\begin{proof}[Proof of Remark \ref{rk:18}]
For the first formula, note that, for any  fixed vector 
$a$ which belongs to $ \R^{n+2}$ when $K=1$  and to $\in \R^{n+2}_1$ 
when $K=-1,$ its orthogonal 
projection on $M_K^{n+1}$ is $a-K\langle a,x\rangle x$, which satisfies:
$$ div(a-K\langle a,x\rangle x) = -(n+1)K \langle x,a\rangle~. $$
Integrating this equality over the interior of an oriented hypersurface
$\Sigma$ yields:
$$ (n+1)K\int_{int(\Sigma)} x dv = - \int_\Sigma Ndv~, $$
where $N$ is the unit exterior normal to $\Sigma$. The formula follows by
taking the limit when $\Sigma$ converges to the boundary of a polyhedron in
$M_K^{n+1}$.

For the second formula, choose $p\in \{ 0,\cdots, n\}$ and
take the limit, as $\Sigma$ converges to the boundary
of $P$, of the Minkowski formula:
$$ \int_\Sigma (n-p)KS_{p}x + (p+1)S_{p+1}N dv = 0~. $$
We include a proof of the Minkowski formulae for the reader's  convenience. Let
$a$ be any fixed vector in $ \R^{n+2}$ in case $K=1$ and in $ \R^{n+2}_1$ in case $K=-1.$ 
Call $Z:=a-K\langle a,x\rangle x - \langle a,N\rangle N$ the projection of $a$
on  $\Sigma$. Finally 
take a local orthonormal frame $\{e_1,\ldots,e_n\}$ on $\Sigma$. Since
$div(T_p)=0$, we have: 
\begin{eqnarray*}
div(T_p(Z))& = & \sum_{i=1}^{n} \langle \nabla_{e_i}Z, T_p{e_i}\rangle\\
& = & \sum_{i=1}^{n} \langle D_{e_i}Z, T_p{e_i}\rangle\\
& = & \sum_{i=1}^{n} - K
\langle a,x\rangle \langle T_p(e_i), e_i\rangle - \langle
a,N\rangle \langle \overline{\nabla}_{e_i}N, T_p(e_i)\rangle\\
& = & - (n-p) KS_p \langle a,x\rangle - (p+1)S_{p+1} \langle a,N\rangle ~,
\end{eqnarray*}
where we used properties (2) and (3) in Lemma 2.2. 
Integrating over $\Sigma$ and replacing $p$ by $p-1$ yields the result.
\end{proof}

\subsection{The Euclidean formulas.}

The proof of the Euclidean version of our vector-valued Schl\"afli formula
follows by an elementary scaling argument from Theorem \ref{tm:general}.
 
\begin{proof}[Proof of Theorem \ref{tm:general-e}]
Consider an Euclidean polyhedron $P$, along with a first-order deformation. 
For each $t\in \R, t>0$, we consider the spherical polyhedron $P_t$ which is
the inverse image of the scaled polyhedron
$tP$ by the projective map:
\begin{eqnarray*}
  S_+^{n+1}:= \{ (x_1, \ldots,x_ {n+1})  \in S^{n+1}  | x_{n+1} > 0 \}
          & \rightarrow & \R^{n+1} \\
  x & \mapsto & (x_1, x_2, \cdots, x_{n+1})/x_{n+2}~.
\end{eqnarray*}
Then $P_t\rightarrow N$ as $t\rightarrow 0$, where $N$ is the ``north pole'' 
of $S^{n+1}$, the point of coordinates $(0,0,\cdots, 0,1)$. 

Let $H$ be a codimension $p$ face of $P$, and let $H_t$ be the 
corresponding faces of $P_t$. As $t\rightarrow 0$, $V(H^*_t)$ tends to a  
constant, while $\pi'(H_t)$ behaves as $V(H_t)$, i.e., as $t^{n+1-p}$.
If $F$ is a codimension $p+2$ face of $P$ and $F_t$ is the corresponding face
of $P_t$, then $V'(F_t^*)$ tends to a constant, while $\pi(F_t)$ scales as
$t^{n-1-p}$; but $\pi(F_t)$ is almost parallel to $N$, so that 
the component of $\pi(F_t)$ orthogonal to $N$ scales as $t^{n-p}$. Finally,
if $G$ is a codimension $p+1$ face of $P$ and $G_t$ is the
corresponding face of $P_t$, then $V(G_t)$ scales as $t^{n-p}$,
while $\nu_{G_t}^*(\pi(G_t^*))$ tends to a constant vector, which is
easily seen to be $\nu_G^*(\pi(G^*))$. So the equation in the 
statement of Theorem \ref{tm:general-e} follows from equation $(E'_p)$
by removing the negligible terms as $t\rightarrow 0$.
\end{proof}

The proof of Theorem \ref{tm:codim1-e} follows by a very similar scaling
argument from Theorem \ref{tm:codim1}, we leave the details to the reader.

\begin{remark}\label{rk:eucl} In the same way, Corollary \ref{cr:quatre} implies that for 
any first order deformation of a  polyhedron $P\subset \R^{n+1}$ and any $p\in\{2,\ldots,n-1\}$:
$$ - \sum_{H\in F_{p+1}} V'(H) \pi(H^*) + \sum_{L\in F_{p+2}} 
V'(L^*) \pi(L)  = 0~. $$
\end{remark}
\section{Non-convex polyhedra}\label{sec:nonconvex}

In this section we show how our formulas extend to non-convex polyhedra. We will do this with 
some details  for the formulas appearing  in Theorem \ref{tm:quatre} 
and Corollary \ref{cr:quatre}. The other formulas, like for instance formula $(E_p)$ in Theorem \ref{tm:general}  can be extended 
in a similar way.  
Since the notion of a dual polyhedron makes sense only for 
convex polyhedra, some care is needed 
in defining the quantities involved in those formulas in the non-convex case.

First, we define a non-convex polyhedron as the union of a finite number of
convex polyhedra $P_1, \cdots, P_N$, with disjoint interiors, but such
that $P_i$ shares a codimension 1 face with $P_{i+1}$ for 
$i\in \{ 1,\cdots, n-1\}$.

Let now $F$ be a face of a convex polyhedron $P.$ 
For reasons that will become clear, we change  here our notation
and denote  its dual face, in the dual polyhedron $P^*,$ by $F_P^*$ 
in order to emphasize the dependence on the polyhedron $P.$

Consider a convex polyhedron $P$ which can be cut into
two convex pieces $P', P'',$ that is   $P=P'\cup P''.$   Let $F$
be a codimension $p$ face of $P, P'$ and $P''$. Note that $P'\cap
P''$ is a convex polyhedron in a hyperplane $\Pi$ and $F^*_{P'\cap P''}$ 
will refer to the dual of the face $F$ 
in $P'\cap P''$ as a polyhedron in $\Pi.$

\begin{prop}\label{nonconvex} We have:
$$ V(F^*_{P'}) + V(F^*_{P''}) = V(F^*_{P}) +
{\frac{V(S^{p-1})}{V(S^{p-2})}} V(F^*_{P'\cap P''})$$
$$\pi(F^*_{P'}) + \pi(F^*_{P''}) = \pi(F^*_{P}) +
\frac{V(S^{p})}{V(S^{p-1})} \pi(F^*_{P'\cap P''}),$$
where $V(S^{k})$ is the volume of the canonical $k$-sphere.
\end{prop}

\begin{proof}
In this setting, $P'\cap P''$ can be considered either as a
 polyhedron with non-empty interior in the
hyperplane $\Pi$ of $M^{n+1}_K$ which contains it or as a 
(degenerate) polyhedron in the whole $(n+1)$-dimensional ambient
space $M^{n+1}_K.$  Denote by $\overline{F}^*_{P'\cap P''}$
the dual face of $F$ in the latter case. Note that this does not coincide with
its dual $F^*_{P'\cap P''}$ in the former case. 

Recall that given a face $F$ of a convex polyhedron $P$  in $S^{n+1}$ (resp. in $H^{n+1}$), 
then $F^*_P$  is the set of points $x\in S^{n+1}$ (resp. in de Sitter space $S^{n+1}_1$) such
that the oriented plane orthogonal to $x$ is a support plane of $P$ at $F$. It follows that:
\begin{eqnarray*}
F^*_{P'} \cap F^*_{P''} &=  &F^*_{P'\cup P''}\\
F^*_{P'} \cup F^*_{P''} &=  &\overline{F}^*_{P'\cap P''}.
\end{eqnarray*}

Consequently:
\begin{eqnarray*} V(F^*_{P'}) + V(F^*_{P''}) &= &V({F}^*_{P}) +
 V(\overline{F}^*_{P'\cap P''})\\
 \pi(F^*_{P'}) + \pi(F^*_{P''}) &= &\pi({F}^*_{P}) +
 \pi(\overline{F}^*_{P'\cap P''}).
 \end{eqnarray*} 
 Now observe that  ${F}^*_{P'\cap P''}$ is a polyhedron (with nonempty interior)  in an equatorial sphere $S^{p-2}$ of dimension
 $(p-2)$ in the unit sphere $S^{p-1}$ of the linear space $vect(\overline{F}^*_{P'\cap P''})$ generated by 
 $\overline{F}^*_{P'\cap P''}$ and that $\overline{F}^*_{P'\cap P''}$
 is obtained by joining the corresponding poles (by shortest geodesic arcs) to each point of ${F}^*_{P'\cap P''}.$  It is then easily checked that:
$$ V(\overline{F}^*_{P'\cap P''}) =
\frac{V(S^{p-1})}{V(S^{p-2})}
V(F^*_{ P'\cap P''}) \quad \text{and}\quad \pi(\overline{F}^*_{P'\cap P''})= \frac{V(S^{p})}{V(S^{p-1})}
\pi(F^*_{ P'\cap P''})$$
and the formulas follow.
\end{proof}

The proof would actually work in the same way for a convex
polyhedron cut into more than two convex pieces, with some
additional terms corresponding to the intersections of more
than two of the pieces.

 It is clear from this proposition that, if the 
formulas in Theorem \ref{tm:quatre} and Corollary \ref{cr:quatre}  apply to $P'$ and to $P''$, then they must apply to
$P$.  For instance, in the formula in Theorem \ref{tm:quatre}, the first term  is:
 \begin{eqnarray*}
(n+1-p) K \sum_{F\in F_{p-1}} V'(F) \pi(F^*_{P}) &=&
(n+1-p) K \sum_{F\in F_{p-1} }V'(F)  \pi(F^*_{P'})+(n+1-p) K \sum_{F\in F_{p-1}} V'(F) \pi(F^*_{P''})-\\
&-&(n+1-p) K \frac{V(S^{p-1})}{V(S^{p-2})}\sum_{F\in F_{p-2}} V'(F) \pi({F}^*_{P'\cap P''}),\\
\end{eqnarray*}
where the sums are over $(n+2-p)$-dimensional faces of $P, P', P''$
and $P'\cap P''$ respectively.  The second term is:
 \begin{eqnarray*}
(n+1-p)K \sum_{G\in F_{p}}V(G^*_P)\pi'(G) & =&(n+1-p)K \sum_{G\in F_{p}}V(G^*_{P'})\pi'(G) +(n+1-p)K \sum_{G\in F_{p}}V(G^*_{P''})\pi'(G)-\\
&-&(n+1-p)K \frac{V(S^{p-1})}{V(S^{p-2})}\sum_{G\in F_{p-1}}V(G^*_{P'\cap P''})\pi'(G),
 \end{eqnarray*}
where the sums are now over $(n+1-p)$-dimensional faces of
$P, P', P''$ and $P'\cap P''$ respectively. The third term is:
 \begin{eqnarray*}
p \sum_{H\in F_{p+1}} V(H) \pi'(H^*_{P}) &=&
p \sum_{H\in F_{p+1} }V(H)  \pi'(H^*_{P'})+p \sum_{H\in F_{p+1}} V(H) \pi'(H^*_{P''})-\\
&-&p \frac{V(S^{p+1})}{V(S^{p})}\sum_{H\in F_{p}} V(H) \pi'({H}^*_{P'\cap P''}),\\
\end{eqnarray*}
where the sums are over $(n-p)$-dimensional faces of $P, P', P''$
and $P'\cap P''$ respectively.  Finally the fourth term is:
 \begin{eqnarray*}
p \sum_{L\in F_{p+2}}V'(L^*_P)\pi(L) & =& p \sum_{L\in F_{p+2}}V'(L^*_{P'})\pi(L) +p \sum_{L\in F_{p+2}}V'(L^*_{P''})\pi(L)-\\
&-& p \frac{V(S^{p+1})}{V(S^{p})}\sum_{L\in F_{p+1}}V'(L^*_{P'\cap P''})\pi(L),
 \end{eqnarray*}
where the sums are now over $(n-1-p)$-dimensional faces of
$P, P', P''$ and $P'\cap P''$ respectively. 

Adding the four terms and using  the formula in Theorem  \ref{tm:quatre} for
$P'$ and for $P''$ leaves us with:
$$
(n-(p-1)) K {{V(S^{p-1})}\over {V(S^{p-2})}} \left(
 \sum_{F\in F_{p-2}} V'(F) \pi({F}^*_{P'\cap P''}) +
\sum_{G\in F_{p-1}}V(G^*_{P'\cap P''})\pi'(G)
\right) - \coupeq  -
p  \frac{V(S^{p+1})}{V(S^{p})}\left(\sum_{H\in F_{p}} V(H) \pi'({H}^*_{P'\cap P''}) +\sum_{L\in F_{p+1}}V'(L^*_{P'\cap P''})\pi(L)  \right). 
$$

Now an elementary computation shows that, for any $k\geq 2$:
$$ {V(S^{k+1})\over V(S^{k})} = {k-1 \over k}
 {V(S^{k-1})\over V(S^{k-2})}, $$ 
so that the remaining term is simply:
$$-
{{V(S^{p-1})}\over {V(S^{p-2})}} 
  \left[ (n-(p-1))
 \left(\sum_{F\in F_{p-2}} V'(F) \pi({F}^*_{P'\cap P''}) +
\sum_{G\in F_{p-1}}V(G^*_{P'\cap P''})\pi'(G)\right) + \right.\coupeq \left. +
(p-1)\left(  \sum_{H\in F_{p}} V(H) \pi'({H}^*_{P'\cap P''}) +\sum_{L\in F_{p+1}}V'(L^*_{P'\cap P''})\pi(L) \right)
 \right]. 
$$
So this vanishes because of  the formula in Theorem  \ref{tm:quatre} applied
to $P'\cap P''$.

Now we can prove Theorem  \ref{tm:quatre} for non-convex polyhedra. First,
we have to define the analogues  of  the quantities $V(F^*)$ and $\pi(F^*)$ for a face $F$ of a non-convex polyhedron $P.$
The polyhedron $P$ can be decomposed into a finite number of
convex polyhedra $P_1, \cdots, P_N$, i.e. $P=\cup_i P_i$ and
the $P_i$ have pairwise disjoint interiors. Then we can
apply the formula of Proposition \ref{nonconvex} (actually the analogous
formula for more than two
polyhedra if necessary) and call $V(F^*)$ and $\pi(F^*)$ the results.
$V(F^*)$ and $\pi(F^*)$ are independent of the decomposition $P=\cup_i P_i$
of $P$ into convex polyhedra $P_i$; if $P=\cup_j P'_j$ is
another decomposition into convex pieces, taking the
finer decomposition $P=\cup_{i,j} P_i\cap P'_j$ and
applying Proposition \ref{nonconvex} shows that the values of $V(F^*)$ and $\pi(F^*)$
obtained by $\cup_i P_i, \cup_j P'_j$ and $\cup_{i,j}P_i\cap P'_j$ are identical.

With this definition, it is not too difficult to prove that
the  formulas in Theorem \ref{tm:quatre} and Corollary \ref{cr:quatre} 
apply to deformations of a non-convex polyhedron $P$.
Start by choosing a decomposition of $P$ into convex
pieces $(P_i)_{1\le i\le N}$, and use the corresponding
 formulas for the $P_i$. Then do as in
the proof above for $P=P'\cup P''$ to get rid of the
terms involving the intersection of two or more of the $P_{i}'s$,
and the result follows. 

\begin{remark} Similarily, the correponding formulas make sense  for non-convex polyhedra in the Euclidean case. 
\end{remark}
\section{Applications}

\subsection{Deformations of polygons}

We consider first the simple result on the possible first-order variations
of the lengths and angles of spherical (or hyperbolic) polygons.

\begin{proof}[Proof of Theorem \ref{tm:polygon-s}]
It follows from Theorem \ref{tm:codim1} that, under any first-order
deformation of a spherical convex polygon, the equation in the statement is
satisfied. The formula also applies to non-convex polygons. Indeed, a 
general polygon can be cut up  into convex ones and we can apply the formula to each of them and sum up
the results. The terms corresponding to added vertices sum up to zero 
(since the total angle remains equal to $2\pi$)
and so do the terms corresponding to the added edges as each such edge 
appears twice with opposite orientations. 

Let $\cP_n$ be the space of polygons with $n$ vertices in $S^2$ for which 
not two vertices are the same point, considered
up to global isometries. It is a manifold of dimension $2n-3$.
The functions $l_i$ and $\theta_i, 1\leq i\leq n$, are defined on $\cP_n$.
At each polygon $p\in \cP$, the differentials of the $l_i$ and the $\theta_i$
generate $T^*_p\cP$, because a polygon is uniquely determined by its edge 
lengths and its dihedral angles. So the $dl_i$ and the $d\theta_i$ generate
a vector subspace of dimension $2n-3$. However they satisfy three linear
relations (those in Theorem \ref{tm:polygon-s}, which are obviously
linearly independent as soon as $p$ has at least 3 vertices). This shows that 
those differentials satisfy no other relation.
\end{proof}

\subsection{Recovering the scalar Schl\"afli formulas}

We now show how the previous considerations allow to recover the higher
Schl\"afli formulas \cite{hsf}. We explain this in the spherical case, 
the hyperbolic case is similar, see the next remark. This will follow 
from the following two formulas.

\begin{lemma} For any first order deformation of a convex polyhedron $P\subset
  S^{n+1}\subset \R^{n+2}$ we have:
\begin{eqnarray*}
\forall p\in \{ 1, \cdots, n\}, 
\sum_{F\in F_{p}} V(F^*) V'(F) - p \sum_{G\in F_{p+1}} \langle
 \nu_G(\pi(G)),\pi(G^*)\rangle & = & 0, \\ 
\forall p\in \{ 0, \cdots, n-1\}, 
(n-p) \sum_{G\in F_{p+1}} \langle \nu_G(\pi(G)),\pi(G^*)\rangle + \sum_{H\in
  F_{p+2}} V'(H^*) V(H) & = & 0~. 
\end{eqnarray*}
\end{lemma}

\begin{proof}
Let $a\in \R^{n+2}$ be any point. Apply 
Remark \ref{local:a1} to the polyhedron $\tau_a(P)$ which 
is contained in the unit sphere centered at $a$. This gives:
$$ \sum_{F\in F_p} V(F^*) (\pi'(F) + V'(F)\ a - \nu_F(\pi(F))) - p \sum_{G\in
  F_{p+1}} \int_{\tau_a(G)} 
\langle \overline{\nu}_{\tau_a(G)}(x), \pi(G^*)\rangle x dv = 0. $$
Applying ($F_p$) to $P$ and taking into account the equation:
$$ \int_{\tau_a(G)}
\langle \overline{\nu}_{\tau_a(G)}(x), \pi(G^*)\rangle x dv = \int_{G}
\langle {\nu}_{G}(x), \pi(G^*)\rangle x dv + \langle \nu_G(\pi(G)),\pi(G^*)\rangle\
a~, $$ 
we conclude that
$$\left( \sum_{F\in F_{p}} V(F^*) V'(F) - p \sum_{G\in F_{p+1}} \langle
  \nu_G(\pi(G)),\pi(G^*)\rangle \right)\ a =0.$$ 
We get the first formula.
The second formula is obtained  in the same way using Remark \ref{local:a2}.
\end{proof}

An immediate consequence  is

\begin{cor}[Higher scalar Schl\"afli formulas]
For any first order deformation of a convex polyhedron 
$P\subset S^{n+1}$ we have, for
each $p\in \{1,\ldots,n-1\}$ 
$$  (n-p) \sum_{F\in F_{p}} V(F^*) V'(F) + p \sum_{H\in F_{p+2}} V'(H^*) V(H)
=0.$$ 
\end{cor}

\begin{remark} {\rm We can recover in a similar way the higher scalar 
Schl\"afli formulas in the hyperbolic case. Indeed,
recall that in the hyperboloid  model of the hyperbolic space, 
the latter is the {\it pseudosphere} of imaginary radius $i,$ centered 
at the origin, in the 
Minkowski space $\R^{n+2}_1.$ Then  one may also consider the same 
sphere but now centered 
at any point $a\in \R^{n+2}_1$ and mimic the arguments in the 
spherical case. We omit the details.} 
\end{remark}

\subsection{Deformations of the sphere as a cone-manifold.}

The formulas found for the deformations of spherical polyhedra 
can be applied to deformations of the sphere which let cone 
singularities appear along a stratified subset. We first consider
the two-dimensional case, the singular locus is then a finite set  of points.

\begin{proof}[Proof of Theorem \ref{tm:cone-s-2}]
We can choose a decomposition of the sphere $(S^2, g_0)$
as a union of interiors of
spherical polygons, with a set of vertices which includes the $v_{i}'s$.
Under the deformation   $g_t$ of $g_0$, those polygons are deformed, 
and it is possible to apply to each of them Theorem \ref{tm:polygon-s}.
However the terms corresponding to the first-order variations of the
lengths of the edges cancel out -- each edge appears twice and the
terms are opposite to each other -- so only the terms corresponding
to the variations of the angles remain. This shows that equation 
(\ref{eq:s2}) has to be satisfied.

Conversely, let $v_1, \cdots, v_n$ be points in $S^2$, and let 
$t_1, \cdots, t_n\in \R$, chosen so that $\sum_{i=1}^n t_i v_i=0$. 
Consider a closed simple polygonal line $p$ which has the $v_i$ among its 
vertices. We can actually suppose that its vertices are $v_1, \cdots,
v_p$, with $p\geq n$. Then the complement of $p$ is the disjoint union
of two disks $D_+$ and $D_-$. We can apply Theorem \ref{tm:polygon-s} 
to each of those disks, taking $\theta'_i=-t_i/2$ if $i\leq n$ 
and $\theta'_i=0$ if $n<i\leq p$, and $l'_i=0$.
Theorem \ref{tm:polygon-s} shows that there exist first-order
deformations of $D_+$ and $D_-$ such that the first-order
variation of the total angle at $v_i$ is $-t_i$. This leads to a
first-order deformation of the sphere -- obtained by gluing $D_+$
and $D_-$ isometrically along their boundary -- as needed.
\end{proof}

\begin{proof}[Proof of Theorem \ref{tm:sphere-codim1}]
The argument used in the proof of Theorem \ref{tm:cone-s-2} can also
be used here, based on Theorem \ref{tm:codim1} applied to the top 
dimensional cells of a cell decomposition of $S^{n+1}$. The terms
corresponding to codimension $1$ faces vanish since each appears
twice, once for each codimension $0$ cell containing the codimension
$1$ cell considered, and those two terms vanish.
\end{proof}

\begin{proof}[Proof of Theorem \ref{tm:sphere}]
Applying formula $(K_p)$ to each of all (maximal dimension) cells
of $\cC$, we obtain that: 
$$ 
\sum_{C\in \cC_0} \left(
(n-p+1)\left( \sum_{F\in F_{p-1}} V'(F) \pi(F^*_C) - \sum_{G\in F_p} V'(G_C^*)
 \pi(G)\right) + \right.\coupeq \left. +
p\left( - \sum_{H\in F_{p+1}} V'(H) \pi(H_C^*) + \sum_{L\in F_{p+2}} 
V'(L_C^*) \pi(L) \right) \right) = 0~, 
$$
where $\cC_0$ is the set of codimension $0$ faces of $\cC$, i.e., the set of
the maximal dimension cells of $\cC$. 
The previous formula can be written as 
$$ 
(n-p+1)\left( \sum_{F\in F_{p-1}} V'(F) \sum_{C\supset F} \pi(F^*_C) - 
\sum_{G\in F_p} \sum_{C\supset G} V'(G_C^*)
 \pi(G)\right) + \coupeq +
p\left( - \sum_{H\in F_{p+1}} V'(H) \sum_{C\supset H} \pi(H_C^*) + 
\sum_{L\in F_{p+2}} \sum_{C\supset L} V'(L_C^*) \pi(L) \right) =0~. 
$$
The proof therefore follows by identifying the different terms 
with the quantities introduced before the statement of the theorem.
\end{proof}

\subsection{Cone singularities in hyperbolic polygons.}

The arguments used above to describe first-order variations of 
the 2-dimensional sphere among cone-surfaces can be applied also
for deformations of a hyperbolic polygon among hyperbolic 
polygons containing cone singularities at some points. This is
of course based on the hyperbolic version of Theorem \ref{tm:polygon-s},
applied to each polygon in a cell decomposition of the hyperbolic
polygon considered, with the singular points among the vertices.
The terms corresponding to the interior edges of the cell
decomposition cancel, and the terms at the interior vertices sum
up to the first-order variations of the total angles at those
vertices. This yields Theorem \ref{tm:polygone-h}.

As mentioned in the introduction, it is also possible to state
corresponding results in higher dimensions, for deformations of
hyperbolic polyhedra in which cone singularities appear along
a stratified subset. We leave the details to the interested reader.

\subsection{Isometric deformations}
Let $P$ be a  polyhedron in $\R^{n+1}, H^{n+1}$ or $S^{n+1}.$ A deformation of $P$ is called  isometric
if each codimension 1 face of $P$ remains congruent to itself (through rigid motions) during the deformation. 
For example deforming a polyhedron through ambient rigid motions induces {\it trivial} isometric deformations.
 $P$ is said {\it flexible} if it admits non trivial isometric deformations and {\it rigid} otherwise. 
A result going back to Legendre and Cauchy states that convex polyhedra are rigid. R. Connelly \cite{connelly}
constructed the first example of a flexible polyhedron in $\R^3.$ To our knowledge there are no known examples in higher dimensions. Isometric deformations have infinitesimal versions: a first order deformation of a polyhedron $P$
is called isometric if it keeps the distance on the codimension 1 faces invariant up to the first order. We can deduce from our formulas some constraints on deforming polyhedra  isometrically. Indeed, as a first order isometric deformation 
keeps the volume of all faces (except the interior of the polyhedron) invariant up to the first order, an immediate consequence of
 Theorem \ref{tm:codim1} and Corollary \ref{cr:quatre} (more precisely their non-convex versions, cf. Section \ref{sec:nonconvex}) and 
Remark \ref{rk:eucl} is:
\begin{cor} Let $P$ be a polyhedron in $\R^{n+1}$ ($K=0$), $H^{n+1}$ ($K=-1$) or
$S^{n+1}$ ($K=1$). Under any isometric first order deformation of $P$ we have, for any $p\in \{2,\ldots,n-1\}:$
$$(n-p+1)K\left( - \sum_{G\in F_p} V'(G^*)
 \pi(G)\right) +
p\left( \sum_{L\in F_{p+2}} 
V'(L^*) \pi(L) \right) = 0~, 
$$
and  $$\sum_{F\in F_2} \theta'(F) \pi(F) =0~. $$
\end{cor}
An analoguous formula can be derived using Theorem \ref{tm:quatre}.

\subsection{Possible future extensions.}

It is quite natural to wonder whether the results stated here
on the first-order deformations of the sphere -- or of hyperbolic
polygons or polyhedra -- extend to corresponding deformations of spherical,
Euclidean or hyperbolic manifolds -- for which cone singularities
appear -- or of cone-manifolds. We do not wish to elaborate much
on this here, except to mention that this indeed seems to be
possible, although more elaborate constructions are needed. It
appears in particular necessary to consider a flat bundle over
the manifold which is deformed, which ``contains'' the relations
which have to be satisfied under the deformations. For flat
(cone-)manifolds this is simply the tangent bundle, but for
spherical or hyperbolic $n+1$-dimensional manifolds it is a
bundle of dimension $n+2$, obtained as the pull-back of the
ambient space by the developing map sending the (non-singular
part of) the (cone-)manifold to $S^{n+1}$ or $H^{n+1}$, respectively.
We hope to come back to this question in a subsequent work. 

\bibliographystyle{alpha}
%\bibliography{../outils/biblio}

\begin{thebibliography}{GBKK{\etalchar{+}}07}

\bibitem[Ale05]{alex}
A.~D. Alexandrov.
\newblock {\em Convex polyhedra}.
\newblock Springer Monographs in Mathematics. Springer-Verlag, Berlin, 2005.
\newblock Translated from the 1950 Russian edition by N. S. Dairbekov, S. S.
  Kutateladze and A. B. Sossinsky, With comments and bibliography by V. A.
  Zalgaller and appendices by L. A. Shor and Yu. A. Volkov.

\bibitem[And01]{anderson-L2}
Michael~T. Anderson.
\newblock {$L\sp 2$} curvature and volume renormalization of {AHE} metrics on
  4-manifolds.
\newblock {\em Math. Res. Lett.}, 8(1-2):171--188, 2001.

\bibitem[BC02]{bezdek-connelly2}
K{\'a}roly Bezdek and Robert Connelly.
\newblock Pushing disks apart---the {K}neser-{P}oulsen conjecture in the plane.
\newblock {\em J. Reine Angew. Math.}, 553:221--236, 2002.

\bibitem[BC04]{bezdek-connelly}
K{\'a}roly Bezdek and Robert Connelly.
\newblock The {K}neser-{P}oulsen conjecture for spherical polytopes.
\newblock {\em Discrete Comput. Geom.}, 32(1):101--106, 2004.

\bibitem[BLP05]{BLP}
Michel Boileau, Bernhard Leeb, and Joan Porti.
\newblock Geometrization of 3-dimensional orbifolds.
\newblock {\em Ann. of Math. (2)}, 162(1):195--290, 2005.

\bibitem[Bon98]{bonahon}
Francis Bonahon.
\newblock A {Schl\"afli}-type formula for convex cores of hyperbolic
  3-manifolds.
\newblock {\em J. Differential Geom.}, 50(1):25--58, 1998.

\bibitem[BP01]{boileau-porti}
Michel Boileau and Joan Porti.
\newblock Geometrization of 3-orbifolds of cyclic type.
\newblock {\em Ast\'erisque}, 272:208, 2001.
\newblock Appendix A by Michael Heusener and Joan Porti.

\bibitem[CMS84]{CMS}
Jeff Cheeger, Werner M{\"u}ller, and Robert Schrader.
\newblock On the curvature of piecewise flat spaces.
\newblock {\em Comm. Math. Phys.}, 92(3):405--454, 1984.

\bibitem[Con77]{connelly}
Robert Connelly.
\newblock {A counterexample to the rigidity conjecture for polyhedra}.
\newblock {\em {Inst. Haut. Etud. Sci., Publ. Math.}}, 47:333--338, 1977.

\bibitem[C06]{csikos} Bal\'azs Csik\'os. 
\newblock {A {Schl\"afli}-type formula for polytopes with curved faces and its application to the Kneser-Poulsen conjecture}. 
\newblock {\em Monatsh. Math.},  147(4): 273--292, 2006.



\bibitem[GBKK{\etalchar{+}}07]{GKRS07}
K.~Grosse-Brauckmann, N.~Korevaar, R.~Kusner, J.~Ratzkin, and J.~Sullivan.
\newblock On nondegeneracy and regularity of the classifying map for coplanar
  cmc surfaces.
\newblock In preparation, 2007.

\bibitem[Glu75]{gluck-generique}
Herman Gluck.
\newblock Almost all simply connected closed surfaces are rigid.
\newblock In {\em Geometric topology (Proc. Conf., Park City, Utah, 1974)},
  pages 225--239. Lecture Notes in Math., Vol. 438. Springer, Berlin, 1975.

\bibitem[KKR06]{KKR06}
N.~Korevaar, R.~Kusner, and J.~Ratzkin.
\newblock On the nondegeneracy of constant mean curvature surfaces.
\newblock {\em Geom. Funct. Anal.}, 16(4):891--923, 2006.

\bibitem[KS06]{volume}
Kirill Krasnov and Jean-Marc Schlenker.
\newblock On the renormalized volume of hyperbolic 3-manifolds.
\newblock {math.DG/0607081. To appear, {\it Comm. Math.Phys.}}, 2006.

\bibitem[KW74]{kazdan-warner}
Jerry~L. Kazdan and F.~W. Warner.
\newblock Curvature functions for compact {$2$}-manifolds.
\newblock {\em Ann. of Math. (2)}, 99:14--47, 1974.

\bibitem[Mil94]{milnor-schlafli}
J.~Milnor.
\newblock The {Schl\"afli} differential equality.
\newblock In {\em Collected papers, vol. 1}. Publish or Perish, 1994.

\bibitem[Reg61]{regge}
T.~Regge.
\newblock General relativity without coordinates.
\newblock {\em Nuovo Cimento (10)}, 19:558--571, 1961.

\bibitem[Rei77]{Reilly}
R.~C. Reilly.
\newblock {Applications of the Hessian operator in a Riemannian manifold.}
\newblock {\em Indiana Univ. Math. J.}, 26:459--472, 1977.

\bibitem[RH93]{RH}
Igor Rivin and Craig~D. Hodgson.
\newblock A characterization of compact convex polyhedra in hyperbolic 3-space.
\newblock {\em Invent. Math.}, 111:77--111, 1993.

\bibitem[Ros93]{rosenberg-constant}
Harold Rosenberg.
\newblock Hypersurfaces of constant curvature in space forms.
\newblock {\em Bull. Sci. Math.}, 117(2):211--239, 1993.

\bibitem[RS99]{sem-era}
Igor Rivin and Jean-Marc Schlenker.
\newblock The {Schl\"afli} formula in {Einstein} manifolds with boundary.
\newblock {\em Electronic Research Announcements of the A.M.S.}, 5:18--23,
  1999.

\bibitem[Sch06]{hmcb}
Jean-Marc Schlenker.
\newblock Hyperbolic manifolds with convex boundary.
\newblock {\em Invent. Math.}, 163(1):109--169, 2006.

\bibitem[SS03]{hsf}
Jean-Marc Schlenker and Rabah Souam.
\newblock Higher {S}chl\"afli formulas and applications.
\newblock {\em Compositio Math.}, 135(1):1--24, 2003.

\end{thebibliography}
\newcommand{\etalchar}[1]{$^{#1}$}
\def\cprime{$'$}

\end{document}